\documentclass[final]{siamltex}

\usepackage{mathtools}
\usepackage{latexsym}%,theorem}
\usepackage{verbatim}
\usepackage{cite}
\usepackage[english]{babel}
\usepackage[latin1]{inputenc}
\usepackage{amsmath}%,amsthm}

\usepackage[usenames, dvipsnames]{xcolor}

\usepackage{enumitem}
\usepackage{amsfonts}
\usepackage{amssymb}
\usepackage{amsbsy}
\usepackage{bm}%,palatino}

\let\leq\leqslant
\let\geq\geqslant
\let\emptyset\varnothing

\newcounter{todocounter}

\newcommand{\todoinnum}[3][]%
{\stepcounter{todocounter}\todo[inline,color=#2!40]{{\bf\small \thetodocounter:} \small#3}}

\usepackage[prependcaption,colorinlistoftodos]{todonotes}

\DeclareMathOperator{\im}{im}

\DeclareMathOperator{\cl}{cl}
\DeclareMathOperator{\rint}{ri}
\DeclareMathOperator{\conv}{conv}
\DeclareMathOperator{\cone}{cone}

\DeclareMathOperator{\dom}{dom}
\DeclareMathOperator{\graph}{gr}
\DeclareMathOperator{\gr}{gr}
\DeclareMathOperator{\lin}{lin}
\DeclareMathOperator{\Lin}{Lin}

\let\leq\leqslant
\let\geq\geqslant
\let\emptyset\varnothing

\newcommand{\calB}{\ensuremath{\mathcal{B}}}
\newcommand{\calC}{\ensuremath{\mathcal{C}}}

\newcommand{\calF}{\ensuremath{\mathcal{F}}}

\newcommand{\calK}{\ensuremath{\mathcal{K}}}

\newcommand{\calN}{\ensuremath{\mathcal{N}}}

\newcommand{\calP}{\ensuremath{\mathcal{P}}}
\newcommand{\calQ}{\ensuremath{\mathcal{Q}}}
\newcommand{\calR}{\ensuremath{\mathcal{R}}}
\newcommand{\calS}{\ensuremath{\mathcal{S}}}
\newcommand{\calT}{\ensuremath{\mathcal{T}}}

\newcommand{\calV}{\ensuremath{\mathcal{V}}}
\newcommand{\calW}{\ensuremath{\mathcal{W}}}

\newcommand{\calY}{\ensuremath{\mathcal{Y}}}

\newcommand{\barq}{\ensuremath{\bar{q}}}

\newcommand{\barx}{\ensuremath{\bar{x}}}
\newcommand{\bary}{\ensuremath{\bar{y}}}
\newcommand{\barz}{\ensuremath{\bar{z}}}
\newcommand{\bmat}{\begin{matrix}}
\newcommand{\emat}{\end{matrix}}
\newcommand{\bbm}{\begin{bmatrix}}
\newcommand{\ebm}{\end{bmatrix}}
\newcommand{\bbma}{\begin{bmatrix*}[r]}
\newcommand{\ebma}{\end{bmatrix*}}
\newcommand{\bpm}{\begin{pmatrix}}
\newcommand{\epm}{\end{pmatrix}}
\newcommand{\bpma}{\begin{pmatrix*}[r]}
\newcommand{\epma}{\end{pmatrix*}}
\newcommand{\bvm}{\begin{vmatrix}}
\newcommand{\evm}{\end{vmatrix}}
\newcommand{\bse}{\begin{subequations}}
\newcommand{\ese}{\end{subequations}}
\newcommand{\beq}{\begin{equation}}
\newcommand{\eeq}{\end{equation}}

\newcommand{\ben}{\begin{enumerate}}
\newcommand{\een}{\end{enumerate}}
\newcommand{\beni}{\begin{enumerate}}
\newcommand{\eeni}{\end{enumerate}}
\newcommand{\bena}{\begin{enumerate}}
\newcommand{\eena}{\end{enumerate}}

\newcommand{\bit}{\begin{itemize}}
\newcommand{\eit}{\end{itemize}}
\newcommand{\bthe}{\begin{theorem}}
\newcommand{\ethe}{\end{theorem}}
\newcommand{\blem}{\begin{lemma}}
\newcommand{\elem}{\end{lemma}}
\newcommand{\brem}{\begin{remark}}
\newcommand{\erem}{\end{remark}}
\newcommand{\bprop}{\begin{proposition}}
\newcommand{\eprop}{\end{proposition}}
\newcommand{\bex}{\begin{example}}
\newcommand{\eex}{\end{example}}
\newcommand{\bas}{\begin{assumption}}
\newcommand{\eas}{\end{assumption}}
\newcommand{\bre}{\begin{remark}}
\newcommand{\ere}{\end{remark}}
\newcommand{\bcor}{\begin{corollary}}
\newcommand{\ecor}{\end{corollary}}
\newcommand{\bdfn}{\begin{definition}}
\newcommand{\edfn}{\end{definition}}
\newcommand{\bcon}{\begin{conjecture}}
\newcommand{\econ}{\end{conjecture}}

\newcommand{\inv}{\ensuremath{^{-1}}}
\newcommand{\nonempty}{\ensuremath{\neq\emptyset}}
\newcommand{\pset}[1]{\ensuremath{\{#1\}}}

\newcommand{\zset}{\ensuremath{\pset{0}}}

\newcommand{\set}[2]{\ensuremath{\left\{#1\mid #2\right\}}}
\newcommand{\res}[2]{\ensuremath{#1|_{#2}}}
\newcommand{\abs}[1]{\ensuremath{| #1 |}}

\newcommand{\inn}[2]{\ensuremath{\langle #1 , #2 \rangle}}

\newcommand{\R}{\ensuremath{\mathbb R}}

\newcommand{\N}{\ensuremath{\mathbb N}}

\newcommand{\BP}{\noindent{\bf Proof. }}
\newcommand{\EP}{\hspace*{\fill} $\blacksquare$\bigskip\noindent}
\newcommand{\qand}{\quad\text{ and }\quad}

\renewcommand{\rint}{\mathrm{ri}}

\newcommand{\B}{\mathfrak{B}}
\newcommand{\He}{H_\mathrm{e}}

\newcommand{\Hi}{H_\mathrm{in}}
\newcommand{\Ho}{H_\mathrm{out}}
\newcommand{\Fo}{\calF(H_\mathrm{out})}
\newcommand{\rrm}{\calR_-}
\newcommand{\rrp}{\calR_+}

\newtheorem{example}[theorem]{Example}
\newtheorem{remark}[theorem]{Remark}
\newtheorem{dfn}[theorem]{Definition}
\newtheorem{assumption}[theorem]{Assumption}

\makeatletter
\newsavebox\myboxA
\newsavebox\myboxB
\newlength\mylenA

\newcommand*\widebar[2][0.75]{%
	\sbox{\myboxA}{$\m@th#2$}%
	\setbox\myboxB\null% Phantom box
	\ht\myboxB=\ht\myboxA%
	\dp\myboxB=\dp\myboxA%
	\wd\myboxB=#1\wd\myboxA% Scale phantom
	\sbox\myboxB{$\m@th\overline{\copy\myboxB}$}%  Overlined phantom
	\setlength\mylenA{\the\wd\myboxA}%   calc width diff
	\addtolength\mylenA{-\the\wd\myboxB}%
	\ifdim\wd\myboxB<\wd\myboxA%
	\rlap{\hskip 0.5\mylenA\usebox\myboxB}{\usebox\myboxA}%
	\else
	\hskip -0.5\mylenA\rlap{\usebox\myboxA}{\hskip 0.5\mylenA\usebox\myboxB}%
	\fi}
\makeatother

\renewcommand{\theenumi}{\roman{enumi}}

\title{A geometric approach to convex processes: from reachability to stabilizability }

\author{Jaap Eising \and M. Kanat Camlibel}

\begin{document}
\maketitle

\renewcommand{\thefootnote}{\fnsymbol{footnote}}

\footnotetext{Jaap Eising is with the Department of Mechanical and Aerospace
	Engineering, University of California, San Diego, USA. \texttt{jeising@ucsd.edu}. M. Kanat Camlibel is with the Bernoulli Institute for Mathematics, Computer Science, and Artificial Intelligence, University of Groningen, The Netherlands. {\tt m.k.camlibel@rug.nl}}

\begin{abstract} This paper studies system theoretic properties of the class of difference inclusions of convex processes. We will develop a framework considering eigenvalues and eigenvectors, weakly and strongly invariant cones, and a decomposition of convex processes. This framework will lead to characterizations of reachability, stabilizability and (null-)controllability for nonstrict convex processes in terms of spectral properties. These characterizations generalize all previously known results regarding for instance linear processes and specific classes of nonstrict convex processes. 
\end{abstract}

\section{Introduction}
This paper deals with fundamental system-theoretic properties of difference inclusions of the form 
\begin{equation}\label{eq:diffinc} x_{k+1} \in H(x_k),\quad k\in\N\tag{$\star$}\end{equation}
where $H$ is a convex process, that is, a set-valued map whose graph is a convex cone. System-theoretic properties of systems of the form \eqref{eq:diffinc} and their continuous-time counterparts have received considerable attention ever since the publication of the seminal paper of Aubin et al. \cite{AFO:86}. 

In particular, the case where $H$ is strict (nonempty everywhere) has been extensively studied in the literature. For this case, \cite{AFO:86} provides necessary and sufficient conditions for reachability in continuous-time. For discrete time systems, the paper \cite{Phat:94} provides similar conditions for both reachability and null-controllability. Also, stabilizability has been fully characterized, still under the strictness assumption, by \cite{Phat:96} for discrete-time and by \cite[Thm. 8.10]{Smirnov:02} for continuous-time. The study of stabilizability has been extended to encompass Lyapunov functions in \cite{Goebel:13}.

All these developments heavily rely on a certain duality relation that follows from the strictness assumption. As shown in \cite{Seeger:01}, this duality relation breaks down for nonstrict convex processes in general. Nevertheless, reachability and null-controllability of particular systems of the form \eqref{eq:diffinc} with nonstrict $H$ have been characterized in \cite{HC:07,HC:08,Kaba:15}. Indeed, these papers study linear systems with convex conic constraints. Such systems can be represented by difference inclusions of the form \eqref{eq:diffinc}. The papers \cite{HC:07,HC:08} provide necessary and sufficient conditions for reachability as well as null-controllability under the assumption that the underlying linear system is right-invertible. This assumption is weaker than the strictness assumption. Later, \cite{Kaba:15} has weakened this right-invertibility assumption by explicitly using the formalism of convex processes. Inspired by this line of research on conically constrained linear systems, \cite{ReachNullc:19} has studied general nonstrict convex processes and characterized reachability and null-controllability under certain domain/image conditions that are weaker than the corresponding assumptions \cite{HC:07,HC:08,Kaba:15} make.

%However, the aforementioned results do not yet recover the linear case. 

In this paper we will study system theoretic properties of general nonstrict convex processes. In particular, we are interested in characterizing reachability, (null-)con-trollability, and stabilizability. For linear systems, these properties can be described in terms of invariant subspaces of the state space (see e.g. \cite{Wonham1985,Basile1992,Trentelman:01}). In order to generalize this towards the setting of \eqref{eq:diffinc}, we will develop a framework based on cones that are invariant under convex processes. Together with new results on eigenvalues and eigenvectors of convex processes from \cite{Eising2021a}, this new framework leads to geometric characterizations of the aforementioned properties in terms of the spectral properties of the dual convex processes under certain domain and image conditions of the convex process at hand. These conditions can be easily verified and are readily satisfied in the case of strict convex processes. In addition, we show that all previously known results can be recovered as special cases from the results of this paper. 

%Indeed, we show that invariance is essential in characterizing . A central ingredient of this paper is the main result of , which provides conditions under which a given weakly invariant cone is guaranteed to contain an eigenvalue. Lastly, recall that duality plays a large role in the analysis of convex processes. Therefore, we will link the concepts of invariance and duality, providing conditions under which weak and strong invariance exchange their roles under duality. 
%
%Thanks to this framework we will develop necessary and sufficient conditions for reachability, null-controllability, and stabilizability of nonstrict convex processes under a certain domain condition within the new framework. These characterizations have a few appealing properties. First of all, they take the form of the well-known Hautus test, that is, they are formulated in terms of eigenvalues of the dual convex process. Secondly, the domain condition is shown to be easily verified. Finally, we show that they generalize all previously known results. Aside from these characterizations, the framework that is developed will allow for a number of generalizations. 

The paper is organized as follows: We begin with a formal problem statement in Section~\ref{sec:problem}. To be able to state the main results, certain preliminary notions/results are introduced in Section~\ref{sec:prelims}. This is followed by the main results as well as their relations to the existing ones in Section~\ref{sec:main}. In Section~\ref{sec:invariance}, we introduce the novel invariance based framework. Section~\ref{sec:proofs} contains the proofs of the main results. Finally, the paper closes with conclusions in Section~\ref{sec:conclusion}.

\section{Problem formulation}\label{sec:problem}

Let $H:\R^n\rightrightarrows \R^n$ be a set-valued map. Its 
{\em graph\/} is defined by
\[ \graph H := \{(x,y)\in \R^n\times \R^n \mid y\in H(x)\}.\]
We say that $H$ is \textit{closed}, \textit{convex}, a \textit{process} or a \textit{linear process} if its graph is closed, convex, a cone or a subspace, respectively.  

Consider the difference inclusion
\begin{equation}\label{eq:incl} x_{k+1} \in H(x_k). \end{equation}
A \textit{trajectory} of \eqref{eq:incl} is a sequence $(x_k)_{k\in\mathbb{N}}$ satisfying \eqref{eq:incl} for every $k\geq0$. In what follows, we will introduce several sets associated with \eqref{eq:incl}.

The \textit{behavior} (see e.g. \cite{Willems:91}) is the set of all trajectories:
\[ \mathfrak{B}(H) := \left\lbrace (x_k)_{k\in\mathbb{N}}\mid \eqref{eq:incl} \textrm{ is satisfied for all }k\in\N \right\rbrace . \] 
For an integer $q\geq1$, we define the $q$\textit{-step behavior} as 
\[ \mathfrak{B}_q(H) := \set{ (x_k)_{k=0}^q} {\eqref{eq:incl} \textrm{ is satisfied for all }k\in\pset{0,1,\ldots,q-1}} . \]
The \textit{feasible set} $\mathcal{F}(H)$ is the set of points from which a trajectory emanates:
\beq
\mathcal{F}(H):= \{\xi \mid \exists (x_k)_{k\in\mathbb{N}}\in \mathfrak{B}(H) \textrm{ with } x_0=\xi \}.\label{e:def feas}
\eeq
The \textit{reachable set} $\mathcal{R}(H)$ is the set of points that can be reached from the origin in finite steps:
\beq
\mathcal{R}(H):= \big\lbrace \xi \mid \exists (x_k)_{k=0}^q \in\mathfrak{B}_q(H) \textrm{ s.t. } x_0=0, x_q=\xi \big\rbrace.\label{e:def:reach}
\eeq  
The \textit{stabilizable set} $\mathcal{S}(H)$ is the set of points from which a stable trajectory exists:
\beq
\mathcal{S}(H) := \{\xi \mid \exists (x_k)_{k\in\mathbb{N}}\in \mathfrak{B}(H) \textrm{ with } x_0=\xi, \lim_{k\rightarrow\infty} x_k= 0  \}.\label{e:def:stab}
\eeq
The \textit{exponentially stabilizable set} $\mathcal{S}_{\mathrm{e}}(H)$ is the set of points from which an exponentially stable trajectory exists:
\beq
\!\!\!\calS_{\mathrm{e}}(H)\! := \{\xi \mid \exists (x_k)_{k\in\mathbb{N}}\in \mathfrak{B}(H),\alpha>0,\mu\in [0,1) \text{ s.t. } x_0=\xi\text{ and }| x_k | \leq \alpha\mu^k \abs{\xi} \} \!\!\label{e:def:e-stab}
\eeq%\forall k\in\N 
where $\abs{\cdot}$ denotes the Euclidean norm.

\noindent The \textit{null-controllable set} $\mathcal{N}(H)$ is the set of points that can be steered to the origin in finite steps:
\beq
\mathcal{N}(H):= \big\lbrace \xi \mid \exists (x_k)_{k=0}^q \in\mathfrak{B}_q(H) \textrm{ s.t. } x_0=\xi, x_q=0 \big\rbrace. \eeq 

All the sets defined above inherit algebraic properties of the set-valued map $H$. In particular, they are all convex cones if $H$ is a convex process and subspaces if $H$ is a linear process. However, they do not retain topological properties from the underlying set-valued map in general. Indeed, none of these sets would be necessarily closed even if $H$ is closed.

We say $H$ is \textit{reachable}, \textit{stabilizable},
%\footnote{What we call stabilizability has been named as weak asymptotic stabilizability in \cite{Phat:96} or weak asymptotic stability in \cite{Gajardo:06}.}
\textit{exponentially stabilizable}, \textit{null-controllable} if $\mathcal{F}(H) \subseteq \mathcal{R}(H)$, $\calF(H)\subseteq \calS(H)$, $\calF(H) \subseteq\calS_{\mathrm{e}}(H)$, $\calF(H)\subseteq\calN(H)$, respectively.

Also, we say $H$ is \textit{controllable} if for all $\xi,\eta\in\calF(H)$ there exist $\ell\geq 0$ and $(x_k)_{k\in\N}\in\mathfrak{B}$ such that $x_0=\xi$ and $x_\ell=\eta$. As the origin belongs to $\calF(H)$, we see that $H$ is controllable if and only if it is both reachable and null-controllable.

The problems we study are to find necessary and sufficient conditions for reachability, (exponential) stabilizability, null-controllability, and controllability of convex processes. 

One of the motivations to study convex processes stems from their link to \textit{constrained linear systems}. To elaborate further on this connection, consider the discrete-time linear input/state/output system given by	
\bse\label{e:lin cons}
\begin{align} \label{eq:ocld2}
x_{k+1} &= Ax_k+Bu_k\\
y_k &= Cx_k+Du_k 
\end{align}
where $k\in\N$, $u_k\in\R^m$ is the input, $x_k\in\R^n$ is the state, $y_k\in\R^p$ is the output and the matrices $A,B,C,D$ are of appropriate dimensions. Suppose that the output of this system is constrained by
\beq
y_k\in\calY
\eeq
for all $k\in\N$ where $\calY\subseteq\R^p$ is a convex cone. 
\ese
Now, define the set-valued map $H:\R^n\rightrightarrows\R^n$ by
\beq\label{e:lin ind H}
H(x) := \{Ax+Bu \mid Cx+Du\in \mathcal{Y} \}. 
\eeq
Since $\calY$ is a convex cone, the set-valued map $H$ is a convex process. This shows that we can view the linear constrained system \eqref{e:lin cons} as
a \textit{difference inclusion} of the form \eqref{eq:incl} where $H$ is given by \eqref{e:lin ind H}.

\section{Preliminaries}\label{sec:prelims}
In this section, we will introduce the notational conventions that will be in force throughout the paper as well as the notions that will be employed in the study of reachability and stabilizability. 

\subsection{Convex cones}
Let $\calS,\calT\subseteq\R^q$ be nonempty sets and $\rho\in\R$. We define $\calS+\calT := \{s+t \mid s\in\calS, t\in\calT\}$ and $\rho \calS := \{\rho s \mid s\in\calS\}$. By convention $\calS+\emptyset = \emptyset$ for every $\calS$ and $\rho\emptyset= \emptyset$ for every $\rho$. The relative interior of $\calS$ is denoted by $\rint(\calS)$. We say that $\calS$ is a \textit{cone} if $\rho x\in\calS$ whenever $x\in\calS$ and $\rho\geq 0$. The conic hull of $\calS$ will be denoted by $\cone(\calS)$. We say that a cone is \textit{finitely generated} if it is a conic hull of finitely many vectors. 

Next, we state three auxiliary results that will be used later. Their proofs are rather elementary and will therefore be omitted.

\begin{lemma}\label{lemm:conespacesequiv}
	Let $\calC\subseteq\mathbb{R}^n$ be a convex cone and let $\mathcal{V}\subseteq \mathcal{W}\subseteq\mathbb{R}^n$ be subspaces. Then, $
(\calC\cap \mathcal{W})+\mathcal{V}=(\calC + \mathcal{V})\cap\mathcal{W}$	
and the following statements are equivalent: 
	\begin{enumerate}
		\item\label{lemmitem:capplus} $(\calC\cap \mathcal{W})+\mathcal{V} = \mathcal{W}$.
		\item\label{lemmitem:pluscap} $(\calC + \mathcal{V})\cap\mathcal{W} = \mathcal{W}$.
		\item\label{lemmitem:subplus} $\mathcal{W}\subseteq \calC+\mathcal{V}$.
		\item\label{lemmitem:plusisplus} $\calC+\mathcal{V} = \calC+\mathcal{W}$.		
	\end{enumerate}	
\end{lemma}

\begin{lemma}\label{l:cone min cone is subspace}
Let $\calC,\calK\subseteq\R^n$ be convex cones. The set $\calC-\calK$ is a subspace if and only if $\rint(\calC)\cap\rint(\calK)\nonempty.$
\end{lemma}

%\BP
%We will prove this case by case.
%
%\eqref{lemmitem:capplus} $\Leftrightarrow$\eqref{lemmitem:pluscap}: As $V\subseteq W$, $(C + \mathcal{V})\cap\mathcal{W}= (C\cap \mathcal{W})+(\mathcal{V}\cap\mathcal{W})=(C\cap \mathcal{W})+\mathcal{V}$. 
%
%\eqref{lemmitem:pluscap}$\Rightarrow$\eqref{lemmitem:subplus}: This implication is immediate.
%
%\eqref{lemmitem:subplus}$\Rightarrow$\eqref{lemmitem:plusisplus} : If $\mathcal{W}\subseteq C+\mathcal{V}$, we can add $C$ to both sides. Using the fact that $C$ is a cone we know $C+C=C$ and hence $C+\mathcal{W}\subseteq C+\mathcal{V}$. As $\mathcal{V}\subseteq \mathcal{W}$, the reverse inclusion is immediate, thus proving $C+\mathcal{V} = C+\mathcal{W}$. 
%
%\eqref{lemmitem:plusisplus} $\Rightarrow$\eqref{lemmitem:pluscap}: Using $C+\mathcal{V} = C+\mathcal{W}$, we can intersect both sides with $\mathcal{W}$. Clearly $(C+\mathcal{W})\cap \mathcal{W}=\mathcal{W}$. Thus proving that $(C + \mathcal{V})\cap\mathcal{W} = \mathcal{W}$. \EP

\begin{lemma}\label{lemm:FGisFD}Let $\calC_\ell$ for $\ell\in\N$ be convex cones with $\calC_\ell\subseteq \calC_{\ell+1}$. If $\bigcup_{\ell=0}^\infty\calC_\ell$ is finitely generated then there exists $q\geq 0$ such that $\calC_{q+\ell}=\calC_q$ for all $\ell\geq 0$ and $\bigcup_{\ell=0}^\infty\calC_\ell=\calC_q$. 
\end{lemma}

%\BP Let $\calC=\cone(\calS)$ where $\calS$ is a finite set of vectors. Then, we know that for each element $s\in\calS$ there is a number $q_s$ such that $s\in \calC_{q_s}$. Since the set $\calS$ is finite, there exists $q\geq 0$ such that $q\geq q_s$ for all $s\in\calS$. As $\calC_q$ is a convex cone containing $\calS$, we have $\calC\subseteq\calC_q$. By definition $\calC_q\subseteq\calC$. Thus, we obtain $\calC=\calC_q$. \EP

%A result we will require a few times in the following exposition is stated in the following lemma: 

\subsection{Convex and linear processes}

%Immediate results from the respective definitions (see also e.g. \cite[Lem. 2.1.2]{Aubin:90}), show that the set-valued map $H$ is a convex process if and only if for all $\alpha\geq 0$, $\beta\geq 0$ and $x,y\in\R^n$:
%\begin{equation} \label{eq:convex process svm}
%		\alpha H(x)+\beta H(y)\subseteq H(\alpha x+\beta y)
%\end{equation} 
%In addition, $H$ is a linear process if and only if for all $\alpha$, $\beta\in\R$ and $x,y\in\R^n$:
%\begin{equation} \label{eq:linear process svm}
%	\alpha H(x)+\beta H(y)\subseteq H(\alpha x+\beta y)
%\end{equation} 

Let $H:\R^n\rightrightarrows \R^n$ be a convex process. For all $x,y\in\R^n$ and $\rho> 0$, we have 
\begin{gather}
H(\rho x)= \rho H(x),\\
H(x)+H(y)\subseteq H(x+y),\\
H(x)=H(x)+H(0),\label{eq: H(x) + H(0)= H(x)}
\end{gather}
and $H(x)$ is convex.  

We define its \textit{domain}, \textit{image} and \textit{kernel} by
\begin{align*}
\dom(H) &:= \{x\in \R^n \mid H(x)\neq \emptyset \}, \\
\im(H) 	&:= \{y\in \R^n \mid \exists x\in \mathbb{R}^n \textrm{ s.t. } y\in H(x)\}, \\
\ker(H) &:= \{x\in \R^n \mid 0\in H(x) \},
\end{align*} 
respectively. If $\dom(H)=\R^n$, we say that $H$ is \textit{strict}.

The \textit{inverse} of $H$ is defined by 
$$
H^{-1}(y):=\{x\in \R^n \mid y\in H(x)\}. 
$$
Clearly, $\dom(H^{-1})=\im(H)$, $\im(H\inv)=\dom(H)$, and 
\begin{equation}\label{eq:def of inverse} \graph (H^{-1}) = \begin{bmatrix} 0 & I_n \\ I_n & 0\end{bmatrix}\graph (H). \end{equation}

For $\lambda \in \mathbb{R}$, we define the set-valued map $H-\lambda I$ by $(H-\lambda I)(x):= \{y - \lambda x \mid y\in H(x)\}$. Then, we have 
\[ \graph (H-\lambda I) = \begin{bmatrix} I_n & 0 \\ -\lambda I_n & I_n\end{bmatrix}\graph (H). \]

We denote the image of a set $\calS$ under $H$ by $H(\calS) = \{ y\in\mathbb{R}^n \mid \exists x\in \calS \textrm{ s.t. } y\in H(x)\}$. Powers of $H$ are defined as follows. By convention, $H^0$ is the identity map, that is $H^0(x):=x$ for all $x\in\R^n$. For $q\geq 1$, we define
\[H^{q+1}(x) := H(H^{q}(x)) \quad \forall x\in \mathbb{R}^n.\]

Clearly, $H(0)$, $\dom(H)$, $\ker(H)$, and $\im(H)$ are all convex cones. In addition, the inverse $H^{-1}$, $H^q$, and $H-\lambda I$ are convex processes for all $q\geq 0$ and $\lambda\in\R$.

Similarly, for a linear process $L:\R^n\rightrightarrows \R^n$, we have that $L(0)$, $\dom(L)$, $\ker(L)$, and $\im(L)$ are all subspaces and the set-valued maps $L\inv$, $L-\lambda I$, $L^q$ are all linear processes. Furthermore, for all $x,y\in\R^n$ and nonzero $\rho\in\R$, we have 
\begin{gather}
L(\rho x)= \rho L(x),\\
L(x)+L(y)= L(x+y),\\
L(x)=L(x)+L(0),\label{eq: L(x) + L(0)= L(x)}
\end{gather}
and $L(x)$ is an affine set. 
 
\subsection{Eigenvalues/vectors of convex processes}
A real number $\lambda$ and nonzero vector $\xi\in\R^n$ form an \textit{eigenpair} of $H$ if $\lambda\xi\in H(\xi)$. In this case $\lambda$ is called an \textit{eigenvalue} and $\xi$ is called an \textit{eigenvector} of $H$ corresponding to (the eigenvalue) $\lambda$. 

For each real number $\lambda$, the convex cone $\ker(H-\lambda I)$ contains all eigenvectors corresponding to $\lambda$ and the origin. This set is called the \textit{eigencone} of $H$ corresponding to $\lambda$. This means that $\lambda$ is an eigenvalue of $H$ if and only if $\ker(H-\lambda I)\neq\zset$. Given a set $\calK$, we define the \textit{spectrum of $H$ with respect to $\calK$} by:
\[ \sigma(H,\calK) := \{ \lambda\in\mathbb{R}\mid \exists \xi\in\calK\setminus\zset \textrm{ such that } \lambda\xi\in H(\xi)\}. \] 
and the \textit{spectrum} by $\sigma(H):=\sigma(H,\R^n)$.

\subsection{Dual processes}
For a nonempty set $\calC\subseteq \R^n $, we define the \textit{negative} and \textit{positive polar cone} by
\[\begin{array}{rl}
\calC^- &:= \{ y\in\R^n \mid \langle x,y\rangle\leq 0 \hspace{1em} \forall x\in \calC\}, \\
\calC^+ &:= \{ y\in\R^n \mid \langle x,y\rangle\geq 0 \hspace{1em} \forall x\in \calC\},
\end{array}\] 
respectively. 

Both the negative and positive polar cones of a set are always closed convex cones. In addition, $\calC$, its closure, convex hull and conic hull have the same polar cones. Furthermore, if $\calC$ is a closed convex cone, then $(\calC^-)^-=\calC$. We also point out that for sets $\calC$ and $\calS$:
\begin{equation} (\calC+\calS)^-  =\calC^-\cap \calS^-, \quad (\calC\cap\calS)^- = \cl (\calC^-+\calS^-). \end{equation}

We define the \textit{negative} and \textit{positive dual} processes $H^-$ and $H^+$ of $H$ by:
\bse\label{eq:def of dual} 
	\begin{align}  p\in H^-(q) &\iff \langle p,x\rangle\geq \langle q,y\rangle,\hspace{5 px} \forall x\in \R^n, \hspace{5 px}\forall y\in H(x),  \\[2mm]
	p\in H^+(q) &\iff \langle p,x\rangle\leq \langle q,y\rangle,\hspace{5 px} \forall x\in \R^n, \hspace{5 px}\forall y\in H(x),
	\end{align}
\ese
respectively. The graphs of the dual processes are related to that of $H$ as follows:
\beq \label{eq:def of dual graph}
	\graph (H^-) = \begin{bmatrix} 0 & I_n \\ -I_n & 0 \end{bmatrix} \big( \graph (H)\big)^-\qquad
	\graph (H^+) = \begin{bmatrix} 0 & I_n \\ -I_n & 0 \end{bmatrix} \big( \graph (H)\big)^+.
\eeq

The following lemma collects properties of the dual processes that will be used later.

The following properties of the dual processes follow from the definitions.

\begin{lemma} \label{lemm:dual proc prop}
	Let $H:\mathbb{R}^n\rightrightarrows\mathbb{R}^n$ be a closed convex process. Then, we have
	\begin{enumerate}
		\item\label{lemm:dual proc prop.1} $H(0)=(\dom(H^+))^+=(\dom (H^-))^-$.
		\item\label{lemm:dual proc prop.3} $(H\inv)^-=(H^+)\inv$.
		\item\label{lemm:dual proc prop.4} $\gr(H^-)=-\gr(H^+)$.
		\item\label{lemm:dual proc prop.5} $\gr\big((H^+)^+\big)=\gr\big((H^-)^-\big)=-\gr(H)$.
		\item\label{lemm:dual proc prop.6} $\gr\big((H^-)^+\big)=\gr\big((H^+)^-\big)=\gr(H)$.
		\item\label{lemm:dual proc prop.2} $\big(\dom(H)\big)^- = -H^-(0) = H^+(0)$.

\item $\big(\im(H-\lambda I)\big)^-=\ker(H^--\lambda I)$ for all $\lambda\in\R$.		
	\end{enumerate}
\end{lemma}

\subsection{Minimal and maximal linear process}

For a cone $\calC$ we define $\lin(\calC) := -\calC\cap \calC$ and $\Lin(\calC)=\calC-\calC$. Note that both $\lin(\calC)$ and $\Lin(\calC)$ are subspaces and, in particular $\lin(\calC)$ is the largest subspace contained in $\calC$ whereas $\Lin(\calC)$ is the smallest subspace that contains $\calC$. 

%We will state a lemma that will be useful later: 
%\begin{lemma}\label{lemm:LinC}
%	Let $C\subseteq\mathbb{R}^n$ be a convex cone and let $\mathcal{V}\subseteq\mathbb{R}^n$ be a subspace. The following are equivalent: 
%	\begin{enumerate}
%		\item\label{lemmitem:subspace} $C +\mathcal{V}$ is a subspace. 
%		\item\label{lemmitem:Lin} $C +\mathcal{V}= \mathrm{Lin}(C)+\mathcal{V}$.
%		\item\label{lemmitem:rinonempty} $\rint (C)\cap \mathcal{V} \neq \emptyset$.
%	\end{enumerate}
%\end{lemma}
%
%\BP The proof is given case by case.
%\eqref{lemmitem:subspace}$\Leftrightarrow$\eqref{lemmitem:Lin}: This is immediate as Lin$(C+\mathcal{V})=$Lin$(C)+\mathcal{V}$. 
%
%\eqref{lemmitem:subspace}$\Leftrightarrow$\eqref{lemmitem:rinonempty}: As it is a cone, $C+\mathcal{V}$ is a subspace if and only if $0\in$ ri$(C+\mathcal{V})$ if and only if ri$(C)\cap \mathcal{V}\neq \emptyset$.\EP

Let $H:\mathbb{R}^n\rightrightarrows \mathbb{R}^n$ be a convex process. Associated with $H$, we define two linear processes $L_-$ and $L_+$ by
\begin{equation}\label{eq:def of L- and L+}
\graph(L_-):=\lin\big(\gr(H)\big)\qand \graph(L_+):=\Lin\big(\graph(H)\big).
\end{equation}

Further, $L_-$ and $L_+$ are, respectively, the largest and smallest (with respect to the graph inclusion) linear processes satisfying
\begin{equation}\label{eq:lh-+}
\graph(L_-)\subseteq\graph(H)\subseteq\graph(L_+).
\end{equation}
We call $L_-$ and $L_+$, respectively, the minimal and maximal linear processes associated with $H$. If $H$ is not clear from context, we will specify it by writing $L_-(H)$ and $L_+(H)$.

The domains of the minimal/maximal linear processes are related to that of $H$:
\begin{equation}\label{eq:domains of L- and L+}
\dom(L_-)\subseteq\lin\big(\dom(H)\big)\qand \dom(L_+)=\Lin\big(\dom(H)\big).
\end{equation}
The reverse inclusion for the former does not hold in general. 

Inverses of minimal/maximal linear processes can be characterized in terms of the inverse of $H$:
\begin{equation}\label{eq:Lmin Lmax under inverse}
L_-(H\inv)=L_-\inv(H)\qand L_+(H\inv)=L_+\inv(H).
\end{equation}

A noteworthy property of the minimal linear process, proven in \cite[Lem. 2.1]{Eising2021a}, is stated next.

\begin{lemma}\label{lemm:H and L powers}
Let $H$ be a convex process and let $L$ be any linear process such that $\graph L \subseteq \graph L_-$. For all $x\in\dom(H)$, $y\in\dom(L)$, we have $H(x+y)=H(x)+L(y)$.\end{lemma}
%\BP
%Let $\zeta\in\dom(H)$ and $\eta\in\dom(L)$. Note that
%$$
%H(\zeta)+L(\eta)
%{\subseteq}H(\zeta)+H(\eta)\overset{\eqref{eq:convex process svm}}{\subseteq}H(\zeta+\eta).
%$$
%For the reverse inclusion, first observe that $\eta\in\dom(L)$ implies that $-\eta\in\dom(L)$ as $L$ is a linear process. Then, we have
%$$
%H(\zeta+\eta)+L(-\eta){\subseteq}H(\zeta+\eta)+H(-\eta)\overset{\eqref{eq:convex process svm}}{\subseteq}H(\zeta).
%$$
%This leads to $
%H(\zeta+\eta)\subseteq H(\zeta)-L(-\eta)=H(\zeta)+L(\eta)$
%where the last equality follows since $L$ is a linear process.
%\EP

Note that for a linear process $L$, the positive and negative duals coincide. Therefore, we denote it by $L^\perp:=L^-=L^+$. The minimal and maximal linear processes associated with a convex process enjoy the following properties that immediately follow from the definitions:
\begin{equation}\label{eq:L(H) vs L(H+-)}
L_-(H^-)=L_+(H)^\perp\qand L_+(H^-)=L_-(H)^\perp.
\end{equation}

\section{Main results}\label{sec:main}
In this section we will state the main results whose proofs can be found in Section~\ref{sec:proofs}.

For linear processes, it can easily be verified that both the feasible and reachable sets can be computed in finitely many steps. For later use, we state this fact below and omit its rather elementary proof. 
\begin{lemma}\label{lemm:linear F and R}
	Let $L:\mathbb{R}^n\rightrightarrows\mathbb{R}^n$ be a linear process. Then, $\calF(L)=\dom(L^n)=L^{-n}(\mathbb{R}^n)$ and $\calR(L)= L^n(0)$ are subspaces. 
\end{lemma}

Let $H:\mathbb{R}^n\rightrightarrows \mathbb{R}^n$ be a convex process. In the rest of the paper, we will use the following shorthand notational conventions:
\beq
\rrm :=\calR\big(L_-(H)\big)\qand \rrp :=\calR\big(L_+(H)\big).
\eeq
Both $\rrm$ and $\rrp$ are subspaces that can be computed in finitely many steps as stated in Lemma~\ref{lemm:linear F and R}.

Our main results will rely on the following \textit{domain condition}:
\beq\label{e:dom cond}
\dom (H)+\rrm \text{ is a subspace}\qand  \dom (H)+\rrm=\dom (H)+\rrp.\tag{DC}
\eeq
Note that the domain condition \eqref{e:dom cond} readily holds whenever $H$ is strict, or $H$ is a linear process. 

Since $\dom(H)$ is a convex cone and $\rrm,\rrp$ are subspaces with $\rrm\subseteq\rrp$, it follows from Lemmas~\ref{lemm:conespacesequiv} and \ref{l:cone min cone is subspace} that the domain condition \eqref{e:dom cond} is equivalent to

\[\rint\big(\dom(H)\big)\cap\rrm\nonempty\qand\rrp\subseteq \dom(H)+\rrm.\]

%\begin{gather*}
%L_-:=,\quad \rrm :=\calR\big(L_-(H)\big), \quad \calS_- = \calS(L_-),\\
%L_+:=L_+(H),\quad \rrp := \calR(L_+), \quad \calS_+=\calS(L_+).
%\end{gather*}

Next, we will introduce two convex processes associated with $H$ that capture the behavior of $H$ inside and outside $\rrp$, respectively. 

We define the \textit{inner process} $\Hi: \mathbb{R}^n \rightrightarrows \mathbb{R}^n$ by 
\[ \graph ( \Hi) := \graph (H) \cap ( \rrp\times\rrp) \]
and the \textit{outer process} $\Ho:\R^n\rightrightarrows\R^n$ by
$$
\graph (\Ho):=\Big(\graph ( H )+\big(\rrp\times \rrp\big)\Big)\cap \big(\calV\times \calV\big)
$$
where $\calV\subseteq\R^n$ is a subspace such that $\rrp\oplus\calV=\R^n$.

We will study these processes in detail in Section~\ref{sec:invariance}. For the moment, we mention only the following result that will be needed to state our main results.

\blem\label{l:ho single valued lin proc}
Let $H$ be a convex process satisfying the domain condition \eqref{e:dom cond}. Then, the outer process $\Ho$ is a single valued linear process, $\calF(\Ho)=\big(\calF(H)+\rrp\big)\cap\calV$ is a subspace and $\Ho(\Fo)\subseteq\Fo$. 
\elem 

We denote the restriction of $\Ho$ to $\Fo$ by $\res{\Ho}{\Fo}$.

%\begin{corollary}\label{cor:Lorthreach}
%A linear process $L$ is reachable if and only if $L^\bot$ is reachable.
%\end{corollary}

Now, we are in a position to state the main results of the paper. We begin with reachability.

\begin{theorem}\label{thm:main-r} Let $H$ be a convex process satisfying the domain condition \eqref{e:dom cond}. Then, the following statements are equivalent: 
	\begin{enumerate}[label=(\roman*)]
		\item\label{thm:main-r1} $H$ is reachable.
		\item\label{thm:main-r2} Both $\Hi$ and $\Ho$ are reachable.
		\item\label{thm:main-r3} All eigenvectors of $\Hi^-$ corresponding to eigenvalues in $[0,\infty)$ belong to $\rrp^\bot$ and $\Fo=\zset$. 
	\end{enumerate}
	Moreover, if $H$ is reachable, then $\calR(H)=\rrp$ and $\calR(H)$ is finitely determined.
\end{theorem}

Theorem~\ref{thm:main-r} captures all existing reachability results for convex processes in the literature as special cases. 

The well-known reachability characterization for \textit{strict} convex processes \cite[Thm. 3.1]{Phat:94} (see also \cite[Thm. 0.4]{AFO:86} for the continuous-time counterpart) follows from Theorem~\ref{thm:main-r}. To see this, note that if $H:\mathbb{R}^n\rightrightarrows\mathbb{R}^n$ is strict, that is $\dom(H)=\R^n$, then the domain condition readily holds and $\calF(H)=\R^n$. In view of Lemma~\ref{l:ho single valued lin proc}, this means that $\calF(\Ho)=\calV$. Therefore, Theorem~\ref{thm:main-r} boils down to $H$ is reachable if and only if $\rrp=\R^n$ and $\Hi^-$ does not have any nonnegative eigenvalues. Note that $\rrp=\R^n$ implies $\Hi=H$. Hence, Theorem~\ref{thm:main-r}, for the strict case, states that $\rrp=\R^n$ and $H^-$ does not have any nonnegative eigenvalues. As what \cite{Phat:94} calls the \textit{rank condition} is equivalent to $\rrp=\R^n$, \cite[Thm. 3.1]{Phat:94} is a special case of Theorem~\ref{thm:main-r}. 

Reachability of convex processes of the form \eqref{e:lin ind H} have been studied in  \cite{HC:07} and \cite{Kaba:15}. While \cite[Thm. V.3]{HC:07} assumes that $\im D+C\calT^*=\R^n$, \cite[Thm. 6.3]{Kaba:15} assumes that $\im D+C\calT^*+\calY=\R^n$. Here $\calT^*$ is the so-called strongly reachable subspace associated with the linear system \eqref{e:lin cons}. These assumptions imply $\dom(H)+\rrm=\R^n$ which, in turn, implies the domain condition \eqref{e:dom cond}. Hence both \cite[Thm. V.3]{HC:07} and \cite[Thm. 6.3]{Kaba:15} are special cases of Theorem~\ref{thm:main-r}. In addition, \cite[Thm. 1]{ReachNullc:19} is a special case of Theorem~\ref{thm:main-r} since it works under the stronger domain condition $\dom(H)+\rrm=\R^n$ as well.

Another noteworthy point is that the results \cite[Thm. 3.1]{Phat:94} and \cite[Thm. V.3]{HC:07} require closedness of the convex processes that they deal with whereas closedness is not assumed by Theorem~\ref{thm:main-r}.  

Last but not least, none of the existing results \cite[Thm. 3.1]{Phat:94}, \cite[Thm. V.3]{HC:07}, \cite[Thm. 6.3]{Kaba:15}, and \cite[Thm. 1]{ReachNullc:19} can directly be applied to nonstrict linear processes. In case $H$ is a linear process, we have $H=L_-=L_+$. Together with Lemma~\ref{lemm:linear F and R}, this implies that the domain condition \eqref{e:dom cond} is readily satisfied for linear processes. Moreover, $\Hi$ is reachable whenever $H$ is linear. Therefore, Theorem~\ref{thm:main-r} asserts that a linear process $H$ is reachable if and only if $\calF(\Ho)=\zset$.

The domain condition \eqref{e:dom cond} proves itself useful also in the context of stabilizability. Even though the next result has a very much parallel statement to that of reachability, its proof is substantially more involved as we will see later. This is mainly because of the different nature of the sets $\calR(H)$ and $\calS(H)$. Indeed, as will be discussed in detail in Section~\ref{sec:invariance}, $\calR(H)$ turns out to be a strongly $H$ invariant set whereas $\calS(H)$ is a weakly $H$ invariant set. 

\begin{theorem}\label{thm:main-s} Let $H$ be a convex process satisfying the domain condition \eqref{e:dom cond}. Then, the following statements are equivalent: 
	\begin{enumerate}[label=(\roman*)]
\item\label{thm:main-s1} $H$ is stabilizable.
\item\label{thm:main-s2} $H$ is exponentially stabilizable.
\item\label{thm:main-s3} Both $\Hi$ and $\Ho$ are exponentially stabilizable.
		\item\label{thm:main-s4} All eigenvectors of $\Hi^-$ corresponding to eigenvalues in $[1,\infty)$ belong to $\rrp^\bot$ and all eigenvalues of the linear map $\res{\Ho}{\Fo}$ are in the open unit disc.
	\end{enumerate}
\end{theorem}

To the best of our knowledge, spectral conditions for stabilizability as stated above have appeared only in \cite[Prob. 8.6.4]{Smirnov:02} (see also \cite[Thm. 8.10]{Smirnov:02} for the continuous-time counterpart). Since \cite[Prob. 8.6.4]{Smirnov:02} deals with strict closed convex processes, the domain condition \eqref{e:dom cond} is automatically satisfied. As such, \cite[Prob. 8.6.4]{Smirnov:02} can be recovered as a special case from Theorem~\ref{thm:main-s}. In \cite[Thm. 5.1]{Gajardo:06}, it is shown that a sufficient condition for stabilizability is that the domain of a convex process admits a particular representation via certain types of eigenvectors of the process itself. Since our stabilizability result stated in terms of the eigenvectors of the dual process instead, it is difficult to compare our result with \cite[Thm. 5.1]{Gajardo:06}. Nevertheless, it should be remarked that \cite[Thm. 5.1]{Gajardo:06} could be applied only if $\calF(H)=\dom(H)$ whereas our result does not require this assumption.

The following will illustrate how the results above can be applied to linear systems with convex conic constraints. Furthermore, it is an example of a problem that could not be resolved by previously known methods. 
\begin{example}\label{ex:specific linear system}
	Consider the linear system
	\[ x_{k+1} = \begin{bmatrix} 1 & 1 &0 \\ 0 & 1 &0 \\ 1 & 0& 0\end{bmatrix} x_k + \begin{bmatrix} 0 \\ 1 \\-1\end{bmatrix}u_k. \]
	Let $a,b\in\mathbb{R}$, and assume that the system is constrained to the convex cone:
	\[ \begin{bmatrix} a & b & 0 \end{bmatrix} x_k \geq 0 \textrm{ and } \begin{bmatrix}1 &-1&-1\end{bmatrix}x_k = 0. \] 
	We are interested in checking for which $a,b$ with $b\neq 0$ the system is reachable or stabilizable. In order to apply the results discussed above, we first write the system as a convex process, by taking
	\[ H(x) := \begin{cases} \begin{bmatrix} 1 & 1 &0 \\ 0 & 1 &0 \\ 1 & 0& 0\end{bmatrix} x + \im\begin{bmatrix} 0 \\ 1 \\-1\end{bmatrix} & \textup{ if } \begin{bmatrix} a & b &0 \\ 1&-1&-1\\ -1&1&1 \end{bmatrix} x\geq 0 , \\ \quad\quad\emptyset &\textup{ otherwise}. \end{cases} \]
	We can find that $L_-$ and $L_+$ are given by 
	\begin{align*} L_-(x) &= \begin{cases} \begin{bmatrix} 1 & 1 &0 \\ 0 & 1 &0 \\ 1 & 0& 0\end{bmatrix} x + \im\begin{bmatrix} 0 \\ 1 \\-1\end{bmatrix} & \textup{ if } \begin{bmatrix} a & b &0 \\ 1&-1&-1 \end{bmatrix} x= 0, \\ \quad\quad\emptyset &\textup{ otherwise}. \end{cases} \\
		L_+(x) &= \begin{cases} \begin{bmatrix} 1 & 1 &0 \\ 0 & 1 &0 \\ 1 & 0& 0\end{bmatrix} x + \im\begin{bmatrix} 0 \\ 1 \\-1\end{bmatrix} & \textup{ if } \begin{bmatrix} 1&-1&-1 \end{bmatrix} x= 0, \\ \quad\quad\emptyset &\textup{ otherwise}. \end{cases}	\end{align*}
	It is now straightforward to check that
	\[ \dom (H) +\calR_- = \dom(H)+\calR_+ = \calR_+ = \im \begin{bmatrix} 1& 0 \\0&1\\1&-1\end{bmatrix}.\]
	As such, the domain condition \eqref{e:dom cond} holds. However, $H$ is not linear, and does not satisfy the assumptions of \cite[Thm. 1]{ReachNullc:19}. As such, previously Furthermore, we can conclude that $\calF(H)\subseteq\calR_+$. This means that $H=\Hi$ and $\calF(\Ho)=\{0\}$. Determining the dual of $\Hi$ yields
	\[\Hi^-(x) =H^-(x) = \begin{cases}  \begin{bmatrix} 1 & 0 & 1\\ 1 & 1 & 0 \\ 0&0&0\end{bmatrix}x + \begin{bmatrix} a \\ b\\0 \end{bmatrix} \mathbb{R}_+ +\im\begin{bmatrix}1\\-1\\-1\end{bmatrix} &  \textup{ if } \begin{bmatrix} 0 &1 &-1 \end{bmatrix} x= 0 , \\\emptyset &\textup{otherwise}. \end{cases}	 \]
	We can conclude that $\lambda\xi\in \Hi^-(\xi)$ if and only if either
	\[ \lambda\in\mathbb{R} \textup{ and }\xi\in \calR_+^\bot= \ker\begin{bmatrix} 1 & 0&1\\0 &1&-1 \end{bmatrix}, \]
	or
	\[ \lambda= (1-\tfrac{a}{b}), \begin{bmatrix}b & b&0 \end{bmatrix}\xi\leq 0  \textup{ and } \begin{bmatrix} 0 &1&-1 \end{bmatrix}\xi= 0.\]
	As such, we can use Theorem~\ref{thm:main-r} to conclude that $H$ is reachable if and only if $\tfrac{a}{b}>1$. In a similar fashion, we see by Theorem~\ref{thm:main-s} that $H$ is stabilizable if and only $\tfrac{a}{b}\geq 0$.
\end{example}

Next, we turn our attention to null-controllability. It is a well-known fact that reachability implies null-controllability for discrete-time linear systems which correspond to linear processes in the framework of this paper. Apart from this particular class, however, reachability does not imply null-controllability in general for convex processes. Interestingly, this implication always holds under the domain condition \eqref{e:dom cond}.

\begin{theorem}\label{t:reach imply null and cont}
Let $H$ be a convex process satisfying the domain condition \eqref{e:dom cond}. If $H$ is reachable, then it is null-controllable. In particular, this means that $H$ is reachable if and only if $H$ is controllable.
\end{theorem} 

One may think that the results on reachability and stabilizability can be extended to null-controllability in the same way. However, as illustrated by \cite[Ex. 5]{ReachNullc:19}, the domain condition \eqref{e:dom cond} is not enough to formulate spectral conditions for null-controllability. Nevertheless, it is still possible to give a spectral characterization by assuming the following 
\textit{image condition}:
\beq\label{e:im cond}
H(\rrp)-\big(\calN(L_-(H))\cap\rrp\big)=\rrp.\tag{IC}
\eeq

\begin{theorem}\label{thm:main-n} Let $H$ be a convex process satisfying the domain condition \eqref{e:dom cond} and the image condition \eqref{e:im cond}. Then, the following statements are equivalent: 
	\begin{enumerate}[label=(\roman*)]
\item\label{thm:main-n1} $H$ is null-controllable.

\item\label{thm:main-n2} Both $\Hi$ and $\Ho$ are null-controllable.
		\item\label{thm:main-n3} All eigenvectors of $\Hi^-$ corresponding to eigenvalues in $(0,\infty)$ belong to $\rrp^\bot$ and the linear map $\res{\Ho}{\Fo}$ is nilpotent. 
	\end{enumerate}
\end{theorem} 

Unlike reachability, null-controllability for convex processes has not been extensively studied in the literature. In \cite[Thm. 3.2]{Phat:94}, the authors assume that both $H$ and $H\inv$ are strict. In that case, both the domain condition \eqref{e:dom cond} and the image condition \eqref{e:im cond} are trivially satisfied. As such, \cite[Thm. 3.2]{Phat:94} is a particular case of 
Theorem~\ref{thm:main-n}. Yet another particular case is \cite[Thm. 2]{ReachNullc:19} which works under the stronger domain condition $\dom(H)+\rrm=\rrp=\R^n$ as well as the stronger image condition $\im(H)+\calN(L_-(H))=\R^n$.

\section{Towards the proofs}\label{sec:invariance}

In this section we will introduce the notions and tools that will be used in the proofs of the main results.

\subsection{Strong and weak invariance}
In the rest of the paper the following invariance notions will play a key role.

\begin{dfn}\label{def:invariance}
	Let $H:\mathbb{R}^n\rightrightarrows \mathbb{R}^n$ be a convex process and $\calC\subseteq\mathbb{R}^n$ be a convex cone. We say that $\calC$ is
	\begin{enumerate}[label=(\roman*)]
		\item \emph{weakly} ${H}$ \emph{invariant} if $H(x)\cap \calC \neq \emptyset$ for all $x\in \calC$.
		\item \emph{strongly} ${H}$ \emph{invariant} if $H(x)\subseteq \calC$ for all $x\in \calC$.
	\end{enumerate}
\end{dfn}

From these definitions the following facts immediately follow.
\begin{lemma}\label{lemm:invariance properties.1}
	Let $H:\mathbb{R}^n\rightrightarrows \mathbb{R}^n$ be a convex process. Then, the following statements hold:
	\ben[label=(\roman*)]
	\item\label{lemmitem:invariance properties.weak and strong} A cone $\calW$ is weakly ${H}$ invariant if and only if $\calW\subseteq H^{-1}(\calW)$. A cone $\calS$ is strongly ${H}$ invariant if and only if $H(\calS)\subseteq \calS$.
	\item\label{lemmitem:invariance properties.strong to weak} If $\calS$ is strongly ${H}$ invariant, then it is also weakly invariant if and only if $\calS\subseteq \dom (H)$.
	\een
\end{lemma}
These two notions of invariance enjoy the following property as was proven in \cite[Lemma 4]{ReachNullc:19}.

\begin{lemma}\label{lemm:invariance properties} Let $H$ be a convex process. If $\calW$ and $\calS$ are, respectively, weakly and strongly ${H}$ invariant, then $\calW\cap \calS$ and $\calS-\calW$ are, respectively, weakly and strongly ${H}$ invariant.
\end{lemma}

Next, we will investigate invariance properties of the feasible, reachable, (exponentially) stabilizable and null-controllable sets. 

\begin{lemma}\label{lemm:feasible largest wHi cone} 
	For a convex process $H$, the feasible set $\calF(H)$ is the largest weakly $H$ invariant convex cone. Moreover, $\calF(H)=H^{-1}(\calF(H))$. 
\end{lemma}

\BP 
Clearly, $\calF(H)$ is a convex cone. By definition, we can see that the feasible set of $H$ is weakly $H$ invariant: If a trajectory exists from $x_0$, there also exists one from any corresponding $x_1$. Therefore $x_1\in H(x_0)\cap\calF(H)$ and hence, by Definition~\ref{def:invariance} the set $\calF(H)$ is weakly $H$ invariant. As any weakly $H$ invariant set naturally allows a trajectory, we can see that $\calF(H)$ is the largest weakly $H$ invariant cone. 

For the second part we know by Lemma~\ref{lemm:invariance properties.1} that $\calF(H)\subseteq H^{-1}(\calF(H))$. It thus suffices to prove the reverse. Applying $H^{-1}$ on both sides, we know $H^{-1}(\calF(H))\subseteq H^{-1}\big( H^{-1}(\calF(H))\big)$ and therefore $H^{-1}(\calF(H))$ is weakly $H$ invariant. As $\calF(H)$ is the largest of such, we have proven the statement.
\EP

Any feasible state is contained in the domain of $H^\ell$ for any $\ell$, hence 
\begin{equation}\label{eq:FH in infinite intersection}
\calF(H)\subseteq \bigcap_{\ell\in\N}H^{-\ell}(\R^n)=\bigcap_{\ell\in\N}\dom(H^{\ell}).
\end{equation}

A case where \eqref{eq:FH in infinite intersection} holds as equality is when $\mathcal{F}(H)= \dom ( H^q) $ for some $q$. In this case we say that $\mathcal{F}(H)$ is \textit{finitely determined}.

\begin{lemma} \label{lemm:Ffdiff} 
The feasible set $\mathcal{F}(H)$ is finitely determined if and only if $\dom (H^q) = \dom (H^{q+1}) \text{ for some }q$.
\end{lemma}

\BP As $\mathcal{F}(H)\subseteq \dom (H^{q+1})\subseteq\dom (H^{q}) $ for all $q\geq 0$, necessity is clear. For sufficiency, let $q$ be such that $\dom (H^q)=\dom (H^{q+1})$ and let $x\in\dom (H^q)$. As $\dom (H^q)=\dom (H^{q+1})$, clearly, there exists $y\in H(x)$ such that $y\in \dom (H^q)$. Thus $y\in H(x)\cap \dom (H^q)$, proving that $\dom (H^q)$ is weakly $H$ invariant. As $\mathcal{F}(H)$ is the largest of such sets, we see that $\mathcal{F}(H) =\dom (H^q)$, proving the lemma. \EP

For the reachable set, we can prove analogous results. To do so, first note that
\begin{equation}\label{e:reach set union}
\calR(H) = \bigcup\limits_{q=0}^\infty H^q(0). 
\end{equation}

\begin{lemma}\label{lemm:reach smallest str inv}
	For a convex process $H$, the reachable set $\calR(H)$ is the smallest strongly $H$ invariant convex cone. Moreover, $\calR(H)=H(\calR(H))$.	
\end{lemma}

\BP Clearly, $\calR(H)$ is a convex cone. Let $\xi\in\calR(H)$. From \eqref{e:reach set union}, we see that $\xi\in H^q(0)$ for some $q\geq 0$. Let $\eta\in H(\xi)$. Then, we have $\eta\in H^{q+1}(0)$. This means that $\eta\in\calR(H)$. Therefore, $\calR(H)$ is strongly $H$ invariant. Let $\calR'$ be a strongly $H$ invariant cone. As $0\in\calR'$, we have  $H^\ell(0)\subseteq\calR'$ for all $\ell\geq 0$. Therefore, $\calR\subseteq\calR'$ and hence $\calR(H)$ is the smallest strongly $H$ invariant convex cone. 

For the second part, note that $H(\calR(H))\subseteq \calR(H)$. Applying $H$ on both sides, we see that $H(\calR(H))$ is strongly $H$ invariant. Since $\calR(H)$ is the smallest of such cones, we see that $\calR(H)\subseteq H(\calR(H))$. Therefore, we have $H(\calR(H))=\calR(H)$. 
\EP

%Note that, unlike for \eqref{eq:FH in infinite intersection}, equality holds immediately here. 
We say $\calR(H)$ is \textit{finitely determined} if $\calR(H)= H^q(0)$ for some $q$. Similar to Lemma \ref{lemm:Ffdiff}, we can state the following.

\begin{lemma}\label{lemm:Rfdiff}
The reachable set $\mathcal{R}(H)$ is finitely determined if and only if $H^q(0) = H^{q+1}(0)$ for some $q$. 
\end{lemma}

\BP Since $H^q(0)\subseteq H^{q+1}(0) \subseteq \mathcal{R}(H)$, necessity is clear. As $H(H^q(0))=H^{q+1}(0)$, we see that $H^q(0)$ is a strongly $H$ invariant cone. As $\calR(H)$ is the smallest of such cones, we see that $\calR(H)=H^q(0)$.
\EP

Similar to feasible and reachable sets, stabilizable sets also enjoy certain invariance properties.

\begin{lemma}\label{l:invariance s}
	Let $H$ be a convex process. Then, $\mathcal{S}_e(H)\subseteq\calS(H)\subseteq\mathcal{F}(H)$. In addition, the sets $\mathcal{S}(H)$ and $\calS_{\mathrm{e}}(H)$ are both weakly $H$ invariant convex cones satisfying $\mathcal{S}(H)=H\inv(\mathcal{S}(H))$ and $\calS_{\mathrm{e}}(H)=H\inv(\calS_{\mathrm{e}}(H))$. 
\end{lemma}

\BP The inclusions are immediate from the definitions. Now let $\xi\in\calS(H)$, then there exists a stable trajectory $(x_k)_{k\in\N}\in\B(H)$ such that $x_0=\xi$. Now clearly $x_1\in H(\xi)\cap \calS(H)$ and thus $\calS(H)$ is weakly $H$ invariant. 

Now let $\xi\in\calS_{\mathrm{e}}(H)$. Therefore there exists $\alpha>0, \mu\in [0,1)$ and a trajectory $(x_k)_{k\in\N}\in\calB(H)$ such that $x_0=\xi$ and $\abs{x_k} \leq \alpha \mu^k |\xi| $ for all $k\geq 0$. If $x_1=0$, then we have $x_1\in H(\xi)\cap \calS_{\mathrm{e}}(H)$. If $x_1\neq0$, let $y_k=x_{k+1}$ for all $k\geq 0$. Clearly, $(y_k)_{k\in\N}\in\calB(H)$ and $y_0=x_1$. Note that  
$|y_k| = |x_{k+1} | \leq \alpha \mu^{k+1} |\xi|= \beta \mu^k |y_0|$ for all $k\geq 0$ where $\beta=\alpha\mu|\xi|/|y_0|$. Therefore $y_0=x_1\in H(\xi)\cap \calS_{\mathrm{e}}(H)$ and thus $\calS_{\mathrm{e}}(H)$ is weakly $H$ invariant. 

Due to Lemma~\ref{lemm:invariance properties.1}.\ref{lemmitem:invariance properties.weak and strong}, we already know that $\mathcal{S}(H)\subseteq H\inv(\mathcal{S}(H))$ and $\calS_{\mathrm{e}}(H)\subseteq H\inv(\calS_{\mathrm{e}}(H))$. What remains to be shown are the reverse inclusions. To do so, let first $\eta\in H\inv(\calS(H))$. Then, there must exist $\xi\in\calS(H)$ such that $\eta\in H\inv(\xi)$. Since $\xi\in\calS(H)$, there exists a stable trajectory $(x_k)_{k\in\N}\in\B(H)$ with $x_0=\xi$. Now, define $y_0=\eta$ and $y_k=x_{k-1}$ for $k\geq 1$. Clearly, $(y_k)_{k\in\N}\in\B(H)$ is a stable trajectory. Therefore, $y_0=\eta\in\calS(H)$. Consequently, we can conclude that $H\inv(\calS(H))\subseteq\calS(H)$. The same argument is still valid if one replaces stability by exponential stability. As such, we also have that $H\inv(\calS_{\mathrm{e}}(H))\subseteq\calS_{\mathrm{e}}(H)$.
\EP

For the exponentially stabilizable set, we have the following property. 
\begin{lemma}\label{lemm:expstabinvariance}
	Let $H$ be a convex process. Then, the exponentially stabilizable set $\calS_{\mathrm{e}}(H)$ is strongly $(H-\mu I)^{-1}$ invariant for all $\mu\in [0,1)$.
\end{lemma}

\BP
Let $\mu\in [0,1)$ and $\xi\in\calS_{\mathrm{e}}(H)$. If $\xi\not\in\dom(H-\mu I)\inv$, then we have $\emptyset=(H-\mu I)\inv(\xi)\subseteq\calS_{\mathrm{e}}(H)$. If $\xi\in\dom(H-\mu I)\inv$, there exists $\eta\in(H-\mu I)\inv(\xi)$. If $\eta=0$, then $\eta\in\calS_e(H)$. Now suppose that $\eta\neq 0$. Since $\xi\in\calS_{\mathrm{e}}(H)$, we know that there exists an exponentially stable trajectory $(x_k)_{k\in\N}\in\B(H)$ with $x_0=\xi$. Define $y_0=\eta$ and $y_k=\mu^k\eta+\sum_{\ell=0}^{k-1}\mu^{k-1-\ell}x_\ell$ for $k\geq 1$. It can be verified that $y_{k+1}\in H(y_k)$ for all $k\in\N$. 

Since $(x_k)_{k\in\N}$ is exponentially stable, there exist $\nu\in[0,1)$ and $\beta\geq 1$ such that $\abs{x_k} \leq \beta\nu^k \abs{\xi}$. Without loss of generality $\mu\leq\nu$, and as such: 
\[ \abs{y_k} \leq \mu^k \abs{\eta} + \sum_{\ell=0}^{k-1} \mu^{k-1-\ell} \beta \nu^\ell\abs{\xi} \leq \mu^k \abs{\eta} + \sum_{\ell=0}^{k} \beta \nu^{k-1}\abs{\xi}\leq \mu^k\abs{\eta} +k\nu^{k-1}\beta\abs{\xi}.\]
Now let $\zeta$ such that $1>\zeta>\nu\geq \mu$.  Note that there exists $\hat{\alpha}$ such that:
\[ \hat{\alpha}\zeta^k \geq \left(\beta\frac{\abs{\xi}}{\abs{\eta}}\right) k\nu^{k-1}, \] 
for all $k\geq 0$. This implies that
\[ \abs{y_k} \leq (1+\hat{\alpha})\zeta^k\abs{\eta}.  \]
As such, $y_0=\eta\in\calS_e(H)$. Consequently, $\calS_{\mathrm{e}}(H)$ is strongly $(H-\mu I)\inv$ invariant.
\EP

For the null-controllable set, we can prove the following invariance properties in a similar fashion to Lemma~\ref{l:invariance s}.

\begin{lemma}\label{l:invariance n}
	Let $H$ be a convex process. Then, $\calN(H)\subseteq\mathcal{F}(H)$. In addition, $\mathcal{N}(H)$ is weakly $H$ invariant and $\mathcal{N}(H)=H\inv(\mathcal{N}(H))$.\end{lemma}

In addition to the reachable set of $H$, we will also consider the reachable sets of $L_-$ and $L_+$. One relation between these respective reachable sets, is given in the following lemma, previously proven as \cite[Lemma 5]{ReachNullc:19}.

\blem\label{lemm:RH invariance}
Let $H$ be a convex process and denote $\calR=\calR(H)$ and $\calR_+=\calR(L_+)$. If $\dom(H)-\calR$ is a subspace, then $\calR-\calR=\calR_+$.
\elem

Next, we investigate some consequences of the domain condition \eqref{e:dom cond}. 

\begin{lemma}\label{lemm:dom he dom h}
	Let $H$ be a convex process and let $\calV$ be a strongly $H$ invariant subspace with $\rrm\subseteq\calV$. If $\dom(H)+\rrm=\dom(H)+\calV$ then $\dom(H^k)+\rrm=\dom(H^k)+\calV$ for all $k\geq 1$.
\end{lemma}

\BP We can use Lemma~\ref{lemm:conespacesequiv} to give an equivalent statement of the implication in the lemma as: 
\[\calV\subseteq \dom (H)+\rrm \implies \calV \subseteq \dom (H^k)+\rrm \quad \forall k\geq 1 \]
To prove this, we begin with defining the convex process $\He:\R^n\rightrightarrows\R^n$, by
\[\gr(\He)=\gr(H)+\rrm\times\zset.\]

We will start by proving that for all $k\geq 1$, we have $\dom(\He^k)=\dom(H^k)+\rrm$. Clearly $\rrm\subseteq\dom (\He^k)$ and $\dom (H^k)\subseteq\dom (\He^k)$. Therefore $\dom(H^k)+\rrm\subseteq \dom(\He^k)$. 

It remains to show that the reverse inclusion holds. This will be achieved by induction on $k$. For $k=1$, note that clearly $\dom (\He)\subseteq\dom(H)+\rrm$. For the induction step assume that $\dom(\He^k)\subseteq\dom(H^k)+\rrm$ for some $k\geq 1$. Let $\xi\in\dom(\He^{k+1})$. Therefore, there exists $\zeta\in\dom(\He^k)$ such that $\zeta\in\He(\xi)$. By the induction hypothesis, we see that $\zeta=\zeta_1-\zeta_2$ where $\zeta_1\in\dom(H^k)$ and $-\zeta_2\in\rrm$. Hence, we obtain $\zeta_1\in H(\xi)+\zeta_2$.

As $\rrm$ is the reachable set of $L_-$, we know there exists $\eta\in\rrm$ such that $\zeta_2\in L_-(\eta)$. This yields $\zeta_1\in H(\xi)+L_-(\eta)$. Then, Lemma~\ref{lemm:H and L powers} implies that $\zeta_1\in H(\xi+\eta)$. Since $\zeta_1\in\dom(H^k)$, we can conclude that $\xi+\eta\in\dom(H^{k+1})$ and hence $\xi\in\dom(H^{k+1})+\rrm$. This proves that for all $k\geq 1$, we have $\dom(\He^k)=\dom(H^k)+\rrm$.

To prove the lemma, recall that it suffices to show that $\calV\subseteq\dom(H^k)+\rrm$. By the hypothesis, we have $\calV\subseteq\dom(H)+\rrm=\dom(\He)$. Next, we claim that $\calV$ is strongly $\He$ invariant. To see this, let $x\in\calV$ and $y\in\He(x)$. Then, there must exist $x'\in\dom(H)$ and $z\in\rrm$ such that $(x,y)=(x',y)-(z,0)$ due to the definition of $\He$. As such, we have $y\in H(x+z)\subseteq H(\calV)\subseteq \calV$ where the first inclusion follows from the fact that $x+z\in\calV+\rrm=\calV$ and the second from the fact that $\calV$ is strongly $H$ invariant. Therefore, we have $\He(\calV)\subseteq\calV$, in other words, $\calV$ is strongly $\He$ invariant. Since $\calV\subseteq\dom(\He)$, we can use Lemma~\ref{lemm:invariance properties.1}.\ref{lemmitem:invariance properties.strong to weak} to conclude that $\calV$ is also weakly $\He$ invariant. Therefore, we get $\calV\subseteq\calF(\He)\subseteq\dom(\He^k)=\dom(H^k)+\rrm$ for all $k\geq 1$ where the first inclusion follows from the fact that $\calF(\He)$ is the largest weakly $\He$ invariant cone, the second from \eqref{eq:FH in infinite intersection}, and finally the third from the first part of this proof.\EP

The image of a convex cone under a convex process enjoys the following duality relation as proven in \cite[Prop. 1]{ReachNullc:19} (also see \cite[Thm. 2.5.7]{Aubin:90} which additionally assumes closedness).
%Note that this result will play a central role in the coming proofs of our main theorem. 

\begin{proposition}\label{prop:duality of images}
	Let $H$ be a convex process and $K$ be a convex cone such that
	$K-\dom(H)$ is a subspace. Then, 
	$$
	\big(H(K)\big)^- =(H^-)^{-1}( K^-).
	$$
\end{proposition}%
Strong and weak invariance become dual notions under certain conditions (see \cite[Thm. 3]{ReachNullc:19}).

\bprop\label{prop:strong inv becomes weak in the dual}
Let $H$ be a convex process and $K$ be a convex cone such that $K-\dom(H)$ is a subspace. Then, $K^-$ is weakly $H^-$ invariant if $K$ is strongly $H$ invariant. Conversely, $K$ is strongly $H$ invariant if $K$ is closed and $K^-$ is weakly $H^-$ invariant. 
\eprop

\subsection{Eigenvalues of convex processes} 
For a detailed study of eigenvalues and eigenvectors of convex processes, we refer the reader to \cite{Eising2021a}. Here, we quote the following result (see \cite[Thm. 3.7]{Eising2021a}) that deals with the location of eigenvectors of a convex process that leaves a certain cone weakly invariant.

\begin{proposition}\label{prop: W = linK}
	Let $H:\mathbb{R}^n\rightrightarrows\mathbb{R}^n$ be a closed convex process and $\calK\subseteq\R^n $ be a weakly $H$ invariant closed convex cone such that $H(0)\cap \calK$ is a subspace, $\lin (\calK)$ is weakly $L_-(H)$ invariant and $\lin(\calK) \subseteq (L_-(H)-\lambda I)\lin (\calK) $ for all $\lambda\geq 0$. Then $\calK=\lin(\calK)$ if and only if any eigenvector of $H$ in $K$ corresponding to an eigenvalue $\lambda\geq 0$ belongs to $\lin K$.
\end{proposition}

Next, we relate nonnegative eigenvalues of $H^-$ to the reachable set of $H$. 

\begin{lemma}\label{lemm:HdualR-} 
	Let $H$ be a convex process and $\lambda\geq 0$. Then, $\ker(H^--\lambda I) \subseteq \calR(H)^-$. 
\end{lemma}

\BP For this, let $\xi$ be an eigenvector of $H^-$ corresponding to a nonnegative eigenvalue $\lambda$, that is, $\lambda\xi\in H^-(\xi)$.  Clearly this means that $(\lambda^{j}\xi,\lambda^{j+1}\xi)\in\graph H^-$ for any $j\geq 0$. Now take $\eta\in \calR(H)$, i.e. $\eta \in H^q(0)$ for some $q$. Then there exists a (finite) sequence $(x_k)_{k=0}^q$ such that $x_0=0$, $x_q=\eta$ and $(x_k,x_{k+1})\in\graph H$ for $k=0,\ldots,q-1$. By the definition of the dual process in \eqref{eq:def of dual}, we know:
\[ \langle \lambda^{j+1}\xi,x_k \rangle \geq \langle \lambda^j\xi, x_{k+1} \rangle  \] 
for any $j\geq 0$ and $k=0,\ldots, q -1$. In particular we can conclude that: 
\[ 0= \langle \lambda^{q}\xi,x_0 \rangle \geq \langle \lambda^{q-1}\xi, x_{1} \rangle \geq \cdots \geq \langle \lambda\xi, x_{q-1}\rangle \geq \langle \xi, x_{q} \rangle=\langle \xi,\eta\rangle . \]
This allows us to conclude that $\xi\in\calR(H)^-$.  \EP

\subsection{Inner and outer processes}
In this section, we will study inner and outer processes as defined in Section~\ref{sec:main}. We begin with recalling the definition of the inner process $\Hi: \mathbb{R}^n \rightrightarrows \mathbb{R}^n$: 
\beq\label{e:def hi}
\graph ( \Hi) := \graph (H) \cap ( \rrp\times\rrp).
\eeq

%Also, we define in addition the process $L_{+,\mathrm{in}}: \mathbb{R}^n\rightrightarrows \mathbb{R}^n$ by:
%\[ \graph ( L_{+,\mathrm{in}} ) := \graph (L_+) \cap ( \rrp\times\rrp). \]

The following observations readily follow from \eqref{e:def hi}:
\bse\label{e:hi props}
\begin{align}
\dom (\Hi^k) &= \dom(H^k)\cap \rrp\,\,\,\forall\,k\geq 1 \label{e:dom hi}\\
\calF(\Hi) & = \calF(H) \cap \rrp,\label{e:f hi}\\
\calR(\Hi) & = \calR(H), \label{e:r hi= r h}\\
\calS(\Hi) & = \calS(H) \cap \rrp, \\
\calS_{\mathrm{e}}(\Hi) & = \calS_{\mathrm{e}}(H) \cap \rrp,\\
\calN(\Hi) &=\calN(H)\cap\rrp\label{e:n hi}.  
\end{align}
\ese

The subsequent results play an instrumental role in studying reachability, stabilizability and null-controllability.

\blem\label{l:Hin and W}
Let $H$ be a convex process satisfying the domain condition \eqref{e:dom cond}. Then, $\calF(\Hi)=\dom(\Hi^n)$ and $\calF(\Hi)+\rrm=\rrp$. Moreover, if $\calW$ is a weakly $\Hi$ invariant convex cone, then the following statements are equivalent:
\begin{enumerate}[label=(\roman*)]
\item\label{l:Hin and W.1} $\calF(\Hi)\subseteq\calR(\Hi)-\calW$
\item\label{l:Hin and W.2} $\calR(\Hi)-\calW=\rrp$
\item\label{l:Hin and W.3} $\Hi^-$ has no eigenvector in $\big(\calR(\Hi)-\calW\big)^-\setminus\rrp^\bot$ that corresponds to a nonnegative eigenvalue. 
\end{enumerate}
If, in addition, $\Hi\inv(\calW)\subseteq\calW$ then the statements above are equivalent to:
\begin{enumerate}[resume,label=(\roman*)]
\item\label{l:Hin and W.4} $\calF(\Hi)=\calW$
\end{enumerate}
\elem

\BP To show that $\calF(\Hi)=\dom(\Hi^n)$, we first claim that $\calR(L_-(\Hi))=\rrm$. Since $\gr(\Hi)\subseteq\gr(H)$, we readily have that $\calR(L_-(\Hi))\subseteq\rrm$. For the reverse inclusion, let $\xi\in\rrm$. Then, there exists $q\geq 0$ and $(x_k)_{k=0}^q\in\B_q(L_-(H))$ such that $x_0=0$ and $x_q=\xi$. Note that $x_k\in\rrm$ for all $k\in\pset{0,1,\ldots,q}$. Since $\rrm\subseteq\rrp$ and $\gr(L_-(\Hi))=\gr(L_-(H))\cap(\rrp\times\rrp)$, we see that $(x_k)_{k=0}^q\in\B_q(L_-(\Hi))$ and hence $\xi\in\calR(L_-(\Hi))$. This proves that
\beq\label{e:L- Hi reach = rrm}
\calR(L_-(\Hi))=\rrm.
\eeq
Next, note that $(\dom(H)\cap\rrp)+\rrm=\rrp$ due to $\rrm\subseteq\rrp$, the domain condition \eqref{e:dom cond} and Lemma~\ref{lemm:conespacesequiv}. In view of \eqref{e:dom hi}, we see that
\beq\label{e:hi dom k nice}
\dom(\Hi)+\rrm=\rrp.
\eeq
Now, we claim that $\dom(\Hi^n)=\dom(\Hi^{n+1})$. To see this, let $\eta\in\dom(\Hi^n)$ and $\zeta\in\Hi^n(\eta)$. Since $\gr(\Hi)\subseteq\rrp\times\rrp$, we see that $\zeta\in\rrp$. From \eqref{e:hi dom k nice}, it follows that $\zeta=\zeta_1-\zeta_2$ where $\zeta_1\in\dom(\Hi)$ and $\zeta_2\in\rrm$. As $\rrm\subseteq\Hi^n(0)$ due to \eqref{e:L- Hi reach = rrm} and Lemma~\ref{lemm:linear F and R}, we have that $\zeta_2\in \Hi^n(0)$. Therefore, $\zeta_1=\zeta+\zeta_2\in\Hi^n(\eta)+\Hi^n(0)=\Hi^n(\eta)$. Since $\zeta_1\in\dom(\Hi)$, we can conclude that $\eta\in\dom(\Hi^{n+1})$. This means that $\dom(\Hi^n)\subseteq\dom(\Hi^{n+1})$. As the reverse inclusion is obvious, we obtain $\dom(\Hi^n)=\dom(\Hi^{n+1})$. Now, it follows from Lemma~\ref{lemm:Ffdiff} that $\calF(\Hi)=\dom(\Hi^n)$.   

To show that $\calF(\Hi)+\rrm=\rrp$, note first that the domain condition \eqref{e:dom cond} and Lemma~\ref{lemm:dom he dom h} imply $\dom(H^n)+\rrm=\dom(H^n)+\rrp$. Since $\rrm\subseteq\rrp$, Lemma~\ref{lemm:conespacesequiv}.\ref{lemmitem:capplus} implies that $(\dom(H^n)\cap\rrp)+\rrm=\rrp$. Then, it follows from \eqref{e:dom hi} that $\dom(\Hi^n)+\rrm=\rrp$. Since $\calF(\Hi)=\dom(\Hi^n)$, we see that $\calF(\Hi)+\rrm=\rrp$. 

For the rest, we will prove the implications \ref{l:Hin and W.1} $\Rightarrow$ \ref{l:Hin and W.2}, \ref{l:Hin and W.2} $\Rightarrow$ \ref{l:Hin and W.1}, \ref{l:Hin and W.2} $\Leftrightarrow$ \ref{l:Hin and W.3}, \ref{l:Hin and W.4} $\Rightarrow$ \ref{l:Hin and W.2}, and finally \ref{l:Hin and W.2} $\Rightarrow$ \ref{l:Hin and W.4} under the extra hypothesis $\Hi\inv(\calW)\subseteq\calW$.

{\em \ref{l:Hin and W.1} $\Rightarrow$ \ref{l:Hin and W.2}}: Since $\rrm\subseteq\calR(H)=\calR(\Hi)$, we see that $\calF(\Hi)+\rrm\subseteq \calR(\Hi)-\calW$. Then, we have $\rrp\subseteq\calR(\Hi)-\calW$ as $\calF(\Hi)+\rrm=\rrp$. The reverse inclusion readily holds since $\calW\subseteq\dom(\Hi)\subseteq\rrp$. Therefore, we can conclude that $\calR(\Hi)-\calW=\rrp$.

{\em \ref{l:Hin and W.2} $\Rightarrow$ \ref{l:Hin and W.1}}: This readily follows from the fact that $\calF(\Hi)\subseteq\rrp$.

{\em \ref{l:Hin and W.2} $\Leftrightarrow$ \ref{l:Hin and W.3}}: For this part of the proof, we want to apply Proposition~\ref{prop: W = linK} to the closed convex cone $\calK:=\big(\calR(\Hi)-\calW\big)^-=\calR(\Hi)^-\cap\calW^+$. To do this, we need to show that the following hypotheses are satisfied:
\begin{enumerate}[label=(\alph*)]
\item\label{it:calk is w hi invariant}
$\calK$ is weakly $\Hi^-$ invariant.
\item\label{it:hi-0 int calk subspace}
$\Hi^-(0)\cap\calK$ is a subspace.
\item\label{it:lin calK invariant}
$\lin(\calK)$ is weakly $L_-(\Hi^-)$ invariant.
\item\label{it:lin calK subset ...}
$\lin(\calK)\subseteq (L_-(\Hi^-)-\lambda I)\lin (\calK) $ for all $\lambda\geq 0$.
\end{enumerate}
To verify these hypotheses, note first that we have
\begin{alignat*}{6}
\dom(\Hi)+\rrm&&\subseteq\,&&\dom(\Hi)-\calR(H)&&\,\subseteq\,&\dom(\Hi)+\rrp\\
%\big(\dom(H)\cap\rrp\big)+\rrm&&\subseteq\,&&\dom(\Hi)-\calR(H)&&\,\subseteq\,&\big(\dom(H)\cap\rrp\big)+\rrp\\
%\big(\dom(H)+\rrm\big)\cap\rrp&&\subseteq\,&&\dom(\Hi)-\calR(H)&&\,\subseteq\,&\rrp\\
\rrp&&\subseteq\,&&\dom(\Hi)-\calR(H)&&\,\subseteq\,&\rrp
\end{alignat*}
where the first line follows from $\rrm\subseteq\calR(H)\subseteq\rrp$ and the second from \eqref{e:hi dom k nice} and \eqref{e:dom hi}. Therefore, we see that $$
\dom(\Hi)-\calR(H)=\rrp.
$$
Since $\calW$ is weakly $\Hi$ invariant, we have 
\beq\label{e:w in rrp}
\calW\subseteq\dom(\Hi)\subseteq\rrp.
\eeq
This means that
\beq\label{e:dom hi minus calk}
\dom(\Hi)-\big(\calR(H)-\calW\big)=\rrp.
\eeq
Since $\calR(H)=\calR(\Hi)$ is strongly $\Hi$ invariant and $\calW$ is weakly $\Hi$ invariant, we see from Lemma~\ref{lemm:invariance properties} that $\calR(H)-\calW$ is strongly $\Hi$ invariant. Then, it follows from \eqref{e:dom hi minus calk} and Proposition~\ref{prop:strong inv becomes weak in the dual} that $\calK$ is weakly $\Hi^-$ invariant. This proves \ref{it:calk is w hi invariant}.

Let $G$ be the convex process given by $\gr(G)=\cl\big(\gr(\Hi)\big)$. Note that $\rint\big(\dom(G)\big)=\rint\big(\dom(\Hi)\big)$. Since $G$ is closed and $G^-=\Hi^-$, Lemma~\ref{lemm:dual proc prop}.\ref{lemm:dual proc prop.1} yields that $\Hi^-(0)=\big(\dom(\Hi)\big)^+$. Then, we see from \eqref{e:dom hi minus calk} that $\Hi^-(0)\cap\calK=\rrp^\bot$ is a subspace. This proves \ref{it:hi-0 int calk subspace}.

From Lemma~\ref{lemm:RH invariance}, we know that $\calR(H)-\calR(H)=\rrp$. Since $\calR(\Hi)=\calR(H)$, it follows from \eqref{e:w in rrp} that $\Lin\big(\calR(H)-\calW\big)=\rrp$. This results in 
\beq\label{e:dual Lin r-w=rrp}
\lin(\calK)=\rrp^\bot.
\eeq
From the definition of $\Hi$ in \eqref{e:def hi}, we have that 
\beq\label{e:rrp bot in hi dual}
\gr(\Hi^-)=\gr(H^-)+(\rrp^\bot\times\rrp^\bot).
\eeq
Note that $\calR_+$ is, by definition, strongly $\Hi$ invariant. Furthermore, $\calR_++\dom(\Hi) = \calR_+$. This means that we can use Proposition~\ref{prop:strong inv becomes weak in the dual} to conclude that $\rrp^\bot$ is weakly $L_-(\Hi^-)$ invariant. Therefore \ref{it:lin calK invariant} holds. 
%In particular, we have $\rrp^\bot\subseteq\Hi^-(\xi)$ for all $\xi\in\rrp^\bot$. This implies that $\rrp^\bot\subseteq(\Hi^-)\inv(\rrp^\bot)$. As such, $\rrp^\bot$ is weakly $\Hi^-$ invariant. Together with \eqref{e:dual Lin r-w=rrp}, this proves \ref{it:lin calK invariant}.

Note that $L_-(\Hi^-)(0)=\lin\big(\Hi^-(0)\big)$. Then, we see that $\rrp^\bot\subseteq L_-(\Hi^-)(0)$ from \eqref{e:rrp bot in hi dual}. Therefore, we have
\begin{align*}
\rrp^\bot\subseteq L_-(\Hi^-)(0)=\big(L_-(\Hi^-)-\lambda I\big)(0)\subseteq\big(L_-(\Hi^-)-\lambda I\big)(\rrp^\bot)
\end{align*}
for all $\lambda\geq 0$. Together with \eqref{e:dual Lin r-w=rrp}, this proves \ref{it:lin calK subset ...}.

Since $\Hi^-$ and $\calK$ satisfy the hypotheses \ref{it:calk is w hi invariant}-\ref{it:lin calK subset ...}, Proposition~\ref{prop: W = linK} and \eqref{e:dual Lin r-w=rrp} imply that \ref{l:Hin and W.2} holds if and only if \ref{l:Hin and W.3} holds.

{\em \ref{l:Hin and W.4} $\Rightarrow$ \ref{l:Hin and W.2}}: As $\calF(\Hi)+\rrm=\rrp$, we see that $\rrp\subseteq\calF(\Hi)-\calR(\Hi)=\calW-\calR(\Hi)$. Since the reverse inclusion is evident, we se that $\calR(\Hi)-\calW=\rrp$.

{\em $\Hi\inv(\calW)\subseteq\calW$ and \ref{l:Hin and W.2} $\Rightarrow$ \ref{l:Hin and W.4}}: Note that $\calR(\Hi)-\calW=\cup_{\ell\geq0}\Hi^\ell(0)-\calW$. Since $\calR(\Hi)-\calW=\rrp$, we see that $\cup_{\ell\geq0}\Hi^\ell(0)-\calW$ is finitely generated. Then, it follows from Lemma~\ref{lemm:FGisFD} that $\Hi^q(0)-\calW=\rrp$ for some $q\geq 0$. Let $\xi\in\calF(\Hi)$. Therefore, there exists a trajectory $(x_k)\in\mathfrak{B}(\Hi)$ such that $x_0=\xi$. Clearly, $x_q\in\Hi^q(\xi)\in\rrp$. Therefore, $x_q=\zeta-\eta$ where $\zeta\in\calW$ and $\eta\in\Hi^q(0)$. This means that $\zeta\in\Hi^q(\xi)$. Thus, we see that $\xi\in\Hi^{-q}(\calW)$. Since $\Hi\inv(\calW)\subseteq\calW$, we further see that $\xi\in\calW$. Therefore, we proved $\calF(\Hi)\subseteq\calW$. The reverse inclusion readily holds since $\calW$ is weakly $\Hi$ invariant and $\calF(\Hi)$ is the largest weakly $\Hi$ invariant set. Therefore, we can conclude that $\calF(\Hi)=\calW$.
\EP

Lemma~\ref{l:Hin and W} leads to the following results for $\Hi$. 

\blem\label{l:Hin char}
Let $H$ be a convex process satisfying the domain condition \eqref{e:dom cond}. Then, the following statements hold:
\begin{enumerate}[label=(\roman*)]
\item\label{l:Hin char.r}$\Hi$ is reachable if and only if all eigenvectors of $\Hi^-$ corresponding to eigenvalues in $[0,\infty)$ belong to $\rrp^\bot$.
\item\label{l:Hin char.s} $\Hi$ is exponentially stabilizable if and only if all eigenvectors of $\Hi^-$ corresponding to eigenvalues in $[1,\infty)$ belong to $\rrp^\bot$.
\item\label{l:Hin char.n1} If $\Hi$ is null-controllable, then all eigenvectors of $\Hi^-$ corresponding to eigenvalues in $(0,\infty)$ belong to $\rrp^\bot$.
\item\label{l:Hin char.n2} Suppose that $H$ satisfies in addition the image condition \eqref{e:im cond}. If all eigenvectors of $\Hi^-$ corresponding to eigenvalues in $(0,\infty)$ belong to $\rrp^\bot$, then $\Hi$ is null-controllable. 
\end{enumerate}
\elem

\BP
{\em\ref{l:Hin char.r}}: By applying Lemma~\ref{l:Hin and W} with the choice $\calW=\zset$, we see that $\Hi$ is reachable if and only if $\Hi^-$ has no eigenvector in $\calR(\Hi)^-\setminus\rrp^\bot$ that correspond to a nonnegative eigenvalue. Since all eigenvectors of $\Hi^-$ corresponding to nonnegative eigenvalues necessarily belong to $\calR(\Hi)^-$ due to Lemma~\ref{lemm:HdualR-}, we can conclude that $\Hi$ is reachable if and only if all eigenvectors of $\Hi^-$ corresponding to eigenvalues in $[0,\infty)$ belong to $\rrp^\bot$.

{\em\ref{l:Hin char.s}}: To prove the `only if' part, let $\lambda\geq 1$ and $\xi$ be such that $\lambda\xi\in\Hi^-(\xi)$. Also let $\barx\in\calF(\Hi)$. Since $\Hi$ is exponentially stabilizable, there exists an exponentially stable trajectory $(x_k)_{k\in\mathbb{N}}\in\mathfrak{B}(\Hi)$ with $x_0=\barx$. Note that we have $\inn{\xi}{x_{k+1}}\leq\lambda\inn{\xi}{x_k}$ for all $k\geq 0$. In particular, this means that
$$
\frac{1}{\lambda^k}\inn{\xi}{x_{k}}\leq\inn{\xi}{x_0}
$$
for all $k\geq 0$. By taking the limit as $k$ tends to infinity, we see that $\xi\in\big(\calF(\Hi)\big)^+=\big(\mathcal{S}_{\mathrm{e}}(\Hi)\big)^+$. Together with Lemma~\ref{lemm:HdualR-}, this implies that $\xi\in\big(\calR(\Hi)-\mathcal{S}_{\mathrm{e}}(\Hi)\big)^-$. From Lemma~\ref{l:invariance s}, we know that $\Hi\inv(\calS_{\mathrm{e}}(\Hi))\subseteq \calS_{\mathrm{e}}(\Hi)$. From Lemma~\ref{l:Hin and W} (implication \ref{l:Hin and W.4} $\Rightarrow$ \ref{l:Hin and W.2} for $\calW=\mathcal{S}_{\mathrm{e}}(\Hi)$), we have that $\calR(\Hi)-\mathcal{S}_{\mathrm{e}}(\Hi)=\rrp$. This proves that $\xi\in\rrp^\bot$.

Finally, what remains to be proven is the `if' part. Since $\Hi\inv(\calS_{\mathrm{e}}(\Hi))\subseteq \calS_{\mathrm{e}}(\Hi)$ due to Lemma~\ref{l:invariance s}, the claim would follow from Lemma~\ref{l:Hin and W} (implication \ref{l:Hin and W.3} $\Rightarrow$ \ref{l:Hin and W.4} for $\calW=\mathcal{S}_{\mathrm{e}}(\Hi)$) if all eigenvectors of $\Hi^-$ within $\big(\calR(\Hi)-\mathcal{S}_{\mathrm{e}}(\Hi)\big)^-$ corresponding to nonnegative eigenvalues belong to $\rrp^\bot$. To show this, suppose, on the contrary, that there exist $\lambda\geq 0$ and
\begin{equation}\label{eq:contradiction} \xi\in\big(\calR(\Hi)-\mathcal{S}_{\mathrm{e}}(\Hi)\big)^-\setminus\rrp^\bot\quad \textup{such that} \quad \lambda\xi\in\Hi^-(\xi).
\end{equation} 
Clearly, we have 
\[ \sigma( \Hi^-,(\calR(\Hi)-\mathcal{S}_{\mathrm{e}}(\Hi))^-\setminus \rrp^\bot) \subseteq \sigma( \Hi^-,\calR(\Hi)^-\setminus \rrp^\bot). \] 
Since $\xi$ belongs to the set on the left hand side, both sets are nonempty. In addition, we know from \cite[Thm. 3.6]{Eising2021a} that both sets are closed and bounded above. Therefore, there exists $\hat{\lambda}\in\sigma( \Hi^-,\calR(\Hi)^-\setminus \rrp^\bot)$ such that $\lambda'\leq\hat{\lambda}$ for all $\lambda'\in\sigma( \Hi^-,\calR(\Hi)^-\setminus \rrp^\bot)$. As all eigenvectors of $\Hi^-$ corresponding to eigenvalues in $[1,\infty)$ belong to $\rrp^\bot$, we know that $\hat{\lambda} < 1$. Now take $\mu$ such that $\hat{\lambda}<\mu<1$. From \eqref{e:rrp bot in hi dual}, we see that $\rrp^\bot\subseteq\ker(\Hi^--\mu I)$. Since $\mu>\hat{\lambda}$, we also see that $\ker(\Hi^--\mu I)\subseteq\rrp^\bot$. Hence, we have 
\beq\label{e:ker hi - mu i}
\ker(\Hi^--\mu I)=\rrp^\bot.
\eeq
From Lemma~\ref{lemm:expstabinvariance}, we know that
\begin{equation}\label{e:hi - mu i inv se}
(\Hi-\mu I)^{-1}(\mathcal{S}_{\mathrm{e}}(\Hi))\subseteq \mathcal{S}_{\mathrm{e}}(\Hi).
\end{equation}
Now, we claim that $\dom\big((\Hi-\mu I)^{-1}\big)-\mathcal{S}_{\mathrm{e}}(\Hi)=\rrp$. To see this, let $G$ be the convex process given by $\gr(G)=\cl\big(\gr(\Hi)\big)$. Since $G$ is closed, it follows from Lemma~\ref{lemm:dual proc prop} that $\big(\im(G-\mu I)\big)^-=\ker(G^--\mu I)$. Note that $\big(\im(G-\mu I)\big)^-=\big(\im(\Hi-\mu I)\big)^-$ and $G^-=\Hi^-$. Therefore, we see that $\big(\im(\Hi-\mu I)\big)^-=\ker(\Hi^--\mu I)$. From \eqref{e:ker hi - mu i} and the fact that $\im(\Hi-\mu I)=\dom\big((\Hi-\mu I)^{-1}\big)$, we see that $\dom\big((\Hi-\mu I)^{-1}\big)=\rrp$. Since $\mathcal{S}_{\mathrm{e}}(\Hi)\subseteq\rrp$, we can conclude that $\dom\big((\Hi-\mu I)^{-1}\big)-\mathcal{S}_{\mathrm{e}}(\Hi)=\rrp$. Then, taking the negative polar of \eqref{e:hi - mu i inv se}, and applying Proposition~\ref{prop:duality of images} to it results in 
\[  \calS_{\mathrm{e}}(\Hi)^- \subseteq \big(\big((\Hi-\mu I)^{-1}\big)^-\big)^{-1} (\calS_{\mathrm{e}}(\Hi)^-).\]
Applying Lemma~\ref{lemm:dual proc prop}.\ref{lemm:dual proc prop.2}, and the fact that $\graph(H^+) = -\graph(H^-)$, we obtain
\[\calS_{\mathrm{e}}(\Hi)^+ = -\calS_{\mathrm{e}}(\Hi)^- \subseteq (\Hi-\mu I)^-(-\calS_{\mathrm{e}}(\Hi)^-).\] 
Lastly, we can use the fact that $ (\Hi-\mu I)^-= (\Hi^--\mu I)$ to show that
\begin{equation}\label{eq:calSHi subset Hi-muI}  \calS_{\mathrm{e}}(\Hi)^+ \subseteq (\Hi^--\mu I) \calS_{\mathrm{e}}(\Hi)^+ \end{equation}

Recall that $\mu>\lambda\geq 0$ and $\xi \in (\calR(\Hi)-\calS_{\mathrm{e}}(\Hi))^-\setminus \rrp^\bot$. As $(\calR(\Hi)-\calS_{\mathrm{e}}(\Hi))^-= \calR(\Hi)^-\cap \calS_{\mathrm{e}}(\Hi)^+$, we see from \eqref{eq:calSHi subset Hi-muI} that there exists $\eta\in \calS_{\mathrm{e}}(\Hi)^+$ such that $\xi+\mu\eta \in \Hi^-(\eta)$. Together with $\lambda\xi \in \Hi^-(\xi)$, this yields   
\[ \mu\big(\xi +(\mu-\lambda)\eta\big) \in \Hi^-(\xi+(\mu-\lambda)\eta).\]
Since $\ker(\Hi^--\mu I)$ is a subspace, we see that $-\xi-(\mu-\lambda)\eta\in\ker(\Hi^--\mu I)$. 
By using the fact that  $\xi+ \mu\eta  \in \Hi(\eta)$, we obtain  
\[ -\lambda\xi \in \Hi^-(-\xi) \] 
Then, Lemma~\ref{lemm:HdualR-} implies that  $-\xi\in\calR(\Hi)^-$. Thus, we see that $\xi\in \calR(\Hi)^-\cap\calR(\Hi)^+=\rrp^\bot$. This contradicts with \eqref{eq:contradiction}.

{\em\ref{l:Hin char.n1}}:  Let $\lambda> 0$ and $\xi$ be such that $\lambda\xi\in\Hi^-(\xi)$. Also let $\barx\in\calF(\Hi)$. Since $\Hi$ is null-controllable, there exists a trajectory $(x_k)_{k\in\mathbb{N}}\in\mathfrak{B}(\Hi)$ with $x_0=\barx$ and $x_q=0$ for some $q\geq 0$. Note that we have $\inn{\xi}{x_{k+1}}\leq\lambda\inn{\xi}{x_k}$ for all $k\geq 0$. In particular, this means that $\inn{\xi}{x_k}\leq\lambda^k\inn{\xi}{\barx}$ for all $k\geq 0$. Since $x_q=0$ and $\lambda>0$, we see that $\inn{\xi}{\barx}\geq 0$. As $\Hi$ is null-controllable, we see that $\xi\in\calN(\Hi)^+$. From Lemma~\ref{lemm:HdualR-}, we know that $\xi\in\calR(\Hi)^-$. As such, we have that $\xi\in\big(\calR(\Hi)-\calN(\Hi)\big)^-$. It follows from Lemma~\ref{l:invariance n} that $\Hi\inv(\mathcal{N}(\Hi))\subseteq\mathcal{N}(\Hi)$. By applying Lemma~\ref{l:Hin and W} (implication \ref{l:Hin and W.4} $\Rightarrow$ \ref{l:Hin and W.2} for $\calW=\mathcal{N}(\Hi)$), we see that $\calR(\Hi)-\mathcal{N}(\Hi)=\rrp$. Therefore, we see that $\xi\in\rrp^\bot$.

{\em\ref{l:Hin char.n2}}: Since $\Hi\inv(\calN(\Hi))\subseteq \calN(\Hi)$ due to Lemma~\ref{l:invariance n}, the claim would follow from Lemma~\ref{l:Hin and W} (implication \ref{l:Hin and W.3} $\Rightarrow$ \ref{l:Hin and W.4} for $\calW=\mathcal{S}_{\mathrm{e}}(\Hi)$) if all eigenvectors of $\Hi^-$ within $\big(\calR(\Hi)-\calN(\Hi)\big)^-$ corresponding to nonnegative eigenvalues belong to $\rrp^\bot$. Since all eigenvectors of $\Hi^-$ corresponding to eigenvalues in $(0,\infty)$ already belong to $\rrp^\bot$, it remains to prove that $\xi\in\rrp^\bot$ whenever $0\in\Hi^-(\xi)$ and $\xi\in\big(\calR(\Hi)-\calN(\Hi)\big)^-$. To see this, note first that $0\in\Hi^-(\xi)$ implies that $\xi\in(\im(\Hi))^-$. Since $\calR(\Hi)\subseteq\im(\Hi)$, we see that $\xi\in\big(\im(\Hi)-\calN(\Hi)\big)^-$. Note that $\im(\Hi)=H(\rrp)$ due to \eqref{e:def hi} and $\calN(\Hi)=\calN(H)\cap\rrp$ due to \eqref{e:n hi}. Since $\calN(L_-(H))\subseteq\calN(H)$, we see that the image condition \eqref{e:im cond} implies that $\rrp\subseteq\im(\Hi)-\calN(\Hi)$. Therefore, $\xi\in\rrp^\bot$.
\EP

The next result will be employed for exponential stabilizability.

\begin{lemma}\label{lemm:poly to stab}
	Let $\calP\subseteq\calF(H)$ be a bounded polyhedron. Suppose that there exist $q$ and $\rho$ with $q\geq 1$ and $\rho\in(0,1)$ such that $H^q(x)\cap \rho \calP \neq \emptyset$ for all $x\in \calP$. Then $\calP\subseteq \calS_{\mathrm{e}}(H)$.
\end{lemma} 

\BP Let $\calP = \conv\{x^i \mid i =1,\ldots, r\}$. Since $H^q(x)\cap \rho \calP \neq \emptyset$ for all $x\in \calP$, there exist 
 $x_0^i,\ldots, x_q^i$ such that
\begin{align}
x^i_0 &= x^i,\\
x^i_q&\in \rho \calP,\label{e:xiq}\\
x^i_{k+1} &\in H(x^i_k)\quad \forall\, k=0,\ldots, q-1.
\end{align} 
From \eqref{e:xiq}, we see that $x^i_q= \rho \sum_j a_{ji} x^j$ where $a_{ji}\geq 0$ and $\sum_j a_{ji}=1$. For $k=0,1,\ldots,q$, define
\[
X_k := \begin{bmatrix} x_k^1 &x_k^2 &\cdots& x_k^r \end{bmatrix}\in \mathbb{R}^{n\times r}
\]
Also, define the matrix $A =(a_{ji})\in\mathbb{R}^{r\times r}$. Then, $X_q = \rho X_0A$. 

Let $\xi\in\calP$. If $\xi=0$, then clearly $\xi\in\calS_{\mathrm{e}}(H)$. Suppose that $\xi\neq 0$. Then, $\xi = \sum_i b_i x^i$  where $b_i\geq 0$ and $\sum_i b_i =1$. Alternatively, we can write $\xi= X_0b$. Now, we construct a trajectory $(x_k)_{k\in\N}$ as follows:
\[ x_{qm+\ell} = \rho^m X_\ell A^mb\quad \textrm{ for }\quad  m \in\N \qand \ell=0,\ldots, q-1.\]
Recall that all entries of $A$ and $b$ are nonnegative. Thus, so are the entries of $A^mb$ for all $m\in\N$. This implies that
\[ x_{qm+\ell+1} = \rho^mX_{\ell+1}A^mb \in H ( \rho^mX_{\ell}A^mb)= H(x_{qm+\ell})\]
for any $m\in\N$ and $\ell=0,\ldots, q-2$. Further, it follows from $X_q = \rho X_0A$ that
\[x_{q(m+1)} = \rho^{m+1}X_0A^{m+1}b = \rho^m X_qA^mb \in H ( \rho^mX_{q-1}A^mb)
=H(x_{qm+(q-1)})
.\]
Therefore $(x_{k})_{k\in\N}\in \mathfrak{B}(H)$ with $x_0=\xi$. It remains to prove that this sequence is exponentially stable. For this, let $\mu$ and $\alpha$ be real numbers such that
\[ \mu^q = \rho, \qand \alpha = \frac{1}{\rho\abs{\xi}} \max_{\substack{1\leq i\leq r \\ 0\leq\ell\leq q-1}} \abs{x_\ell^i} . \]
Then, we see that
\[ \abs{x_{qm+\ell}} = \abs{ \rho^m X_\ell A^mb}  \leq \mu^{mq} \max_{\substack{1\leq i\leq r \\ 0\leq\ell\leq q-1}}  \abs{x_\ell^i}  \leq  \alpha \mu^{qm+\ell} \abs{\xi} \] 
since $\rho\leq\mu^\ell$ for all $\ell$ with $0\leq\ell\leq q-1$. This proves that $\xi\in \mathcal{S}_{\mathrm{e}}(H)$ and therefore concludes the proof. \EP

Now, we turn our attention to the outer process. Recall that outer process $\Ho:\R^n\rightrightarrows\R^n$ is defined by
\begin{equation}\label{eq:def H out}
\graph (\Ho):=\Big(\graph ( H )+\big(\rrp\times \rrp\big)\Big)\cap \big(\calV\times \calV\big)
\end{equation}
where $\calV\subseteq\R^n$ is a subspace such that 
\begin{equation}\label{eq:direct sum R(L+) and W}
\rrp\oplus\calV=\R^n.
\end{equation}
 
Even though the subspace $\calV$ is not unique in general, the subsequent theory will work regardless of the choice made. 

In addition to $\Ho$, we define $L_{+,\mathrm{out}}:\R^n\rightrightarrows\R^n$ in a similar fashion:
\begin{equation}\label{eq:def L out}
\graph (L_{+,\mathrm{out}}):=\Big(\graph (L_+)+\big(\rrp\times \rrp\big)\Big)\cap \big(\calV\times \calV\big).
\end{equation}

The following lemma will collect essential properties of the outer process. 

\begin{lemma}\label{lemm:Houtthings}
	Let $H$ be a convex process. Then, the following statements hold:
	\begin{enumerate}[label=(\roman*)]
		\item\label{lemmitem:domHout} $\dom(\Ho)=\big(\dom( H )+\rrp\big)\cap \calV$.
		\item\label{lemmitem:Loutsinglevalued} $L_{+,\mathrm{out}}$ is a single valued linear process, i.e. a linear map.
	\end{enumerate}
	If, in addition, $H$ satisfies \eqref{e:dom cond}, then we have:
	\begin{enumerate}[resume,label=(\roman*)]
		\item\label{lemmitem:ho=L+out} $\Ho=L_{+,\mathrm{out}}$.
		\item\label{lemmitem:ho^k} $\dom(\Ho^k)=\big(\dom(H^k)+\rrp\big)\cap\calV$ for all $k\geq 1$.
		\item\label{lemmitem:feasible H} $\calF(H)=\dom(H^n)$.
		\item\label{lemmitem:F hout} $\calF(\Ho)=\big(\calF(H)+\rrp\big)\cap\calV$ is a subspace.
		\item\label{lemmitem: calfh + rrp is subspace}
$\calF(H)+\rrp$	is a subspace.		
	\end{enumerate} 
\end{lemma}
 
\BP {\ref{lemmitem:domHout}}: Clearly, we have $\dom(\Ho)\subseteq\big(\dom( H )+\rrp\big)\cap \calV$. As such, it suffices to prove the reverse inclusion $\big(\dom( H )+\rrp\big)\cap \calV \subseteq\dom(\Ho)$. Let $\xi\in \big(\dom( H )+\rrp\big)\cap \calV$. Then, $\xi=\xi_1+\xi_2$ where $\xi_1\in\dom( H )$ and $\xi_2\in\rrp$. Let $\zeta_1\in  H (\xi_1)$. From the direct sum \eqref{eq:direct sum R(L+) and W}, we see that $\zeta_1=\zeta-\zeta_2$ for some $\zeta\in\calV$ and $\zeta_2\in\rrp$. Note that $(\xi,\zeta)=(\xi_1,\zeta_1)+(\xi_2,\zeta_2)\in\graph (\Ho)$. This means that $\xi\in\dom(\Ho)$, proving the statement.

{\ref{lemmitem:Loutsinglevalued}}: For a linear process to be single valued, it is enough to show that $L_{+,\mathrm{out}}(0)=\zset$. To do so, let $y\in L_{+,\mathrm{out}}(0)$. This means that $y\in\calV$ and there exist $x'\in\dom(L_+)\cap\rrp$, $y'\in L_+(x')$, and $z'\in\rrp$ such that $(0,y)=(x',y')+(-x',z')$. Since $\rrp$ is strongly invariant, we know $L_+(\rrp)\subseteq\rrp $ and thus we see that $y'\in \rrp$. This means that $y=y'+z'\in\rrp$. Recall that $y\in\calV$ as well, and hence $y=0$. 

{\ref{lemmitem:ho=L+out}}: It follows from \ref{lemmitem:domHout} that $\dom (L_{+,\mathrm{out}})=\big(\dom ( L_+ )+\rrp\big)\cap \calV$ by replacing $H$ and $\Ho$, respectively, by $L_+$ and $L_{+,\mathrm{out}}$. Since $\dom (H)+\rrp$ is a subspace and $\dom ( L_+ )=\Lin(\dom (H))$ due to \eqref{eq:domains of L- and L+}, we know $\dom ( L_+ )+\rrp= \Lin (\dom (H))+\rrp = \dom (H)+\rrp$. As $L_{+,\mathrm{out}}$ is single valued and its graph is larger than that of $\Ho$, we see that $\Ho$ and $L_{+,\mathrm{out}}$ coincide.

{\ref{lemmitem:ho^k}}: First, we claim that $\big(\dom(H^k)+\rrp\big)\cap\calV\subseteq\dom(\Ho^k)$ for all $k\geq 1$. To see this, let $k\geq 1$ and $x_0\in\big(\dom(H^k)+\rrp\big)\cap\calV$. Then, there exist $y_0\in\dom(H^k)$ and $z_0\in\rrp$ such that $x_0=y_0+z_0$. Since $y_0\in\dom(H^k)$, there exist $y_1,y_2,\ldots,y_k$ such that $y_{\ell+1}\in H(y_\ell)$ with $\ell=0,1,\ldots,k-1$. From the direct sum \eqref{eq:direct sum R(L+) and W}, we see that there exist $z_1,z_2,\ldots,z_k\in\rrp$ and $x_1,x_2,\ldots,x_k\in\calV$ such that $y_\ell=x_\ell-z_\ell$ for all $\ell=1,2,\ldots,k$. Note that $(x_\ell,x_{\ell+1})=(y_\ell,y_{\ell+1})+(z_\ell,z_{\ell+1})\in\gr(\Ho)$ for all $\ell=0,1,\ldots,k-1$ since $(y_\ell,y_{\ell+1})\in\gr(H)$ and $(z_\ell,z_{\ell+1})\in\rrp\times\rrp$. Therefore, we have $x_{\ell+1}\in \Ho(x_{\ell})$ for all $\ell=0,1,\ldots,k-1$ and hence $x_k\in \Ho^k(x_0)$. In other words, $x_0\in\dom(\Ho^k)$. Consequently, we obtain $\big(\dom(H^k)+\rrp\big)\cap\calV\subseteq\dom(\Ho^k)$. 

Therefore, it remains to show that $\dom(\Ho^k)\subseteq\big(\dom(H^k)+\rrp\big)\cap\calV$ for all $k\geq 1$. We will prove this by induction on $k$. Note that $\dom(\Ho)=\big(\dom(H)+\rrp\big)\cap\calV$ due to \ref{lemmitem:domHout}. As the induction hypothesis, we assume that $\dom(\Ho^k)=\big(\dom(H^k)+\rrp\big)\cap\calV$ for some $k\geq 1$. 
 
Now let $x\in\dom(\Ho^{k+1})$. In particular, we have $x\in\dom(\Ho)$. Due to \ref{lemmitem:domHout} and the domain condition, we see that $x_0\in\big(\dom(H)+\rrm\big)\cap\calV$. Therefore, there exist $x_1\in\dom(H)$ and $x_2\in\rrm$ such that $x=x_1+x_2$. Let $y_1\in H(x_1)$. In view of the direct sum \eqref{eq:direct sum R(L+) and W}, $y_1=y-y_2$ where $y\in\calV$ and $y_2\in\rrp$. Then, we have $(x,y)=(x_1,y_1)+(x_2,y_2)\in\gr(\Ho)$ since $(x_1,y_1)\in\gr(H)$ and $(x_2,y_2)\in\rrm\times\rrp\subseteq\rrp\times\rrp$. 
 
Recall that $\Ho$ is single valued due to \ref{lemmitem:Loutsinglevalued} and \ref{lemmitem:ho=L+out}. As $x\in\dom(\Ho^{k+1})$, this means that $y\in\dom(\Ho^k)$. From the induction hypothesis, we then have $y\in\big(\dom(H^k)+\rrp\big)\cap\calV$. Note that $y_1=y-y_2\in\dom(H^k)+\rrp$. Since $\dom(H^k)+\rrp=\dom(H^k)+\rrm$ due to Lemma~\ref{lemm:dom he dom h}, there exist $\xi\in\dom(H^k)$ and $\eta\in\rrm$ such that $y_1=\xi-\eta$. As $\rrm=L_-^n(0)$ there exists $\zeta\in\rrm$ such that $\eta\in L_-(\zeta)$. Therefore, we have $\xi\in H(x_1)+\eta\subseteq H(x_1)+L_-(\zeta)=H(x_1+\zeta)$. Since $\xi\in\dom(H^k)$, we get $x_1+\zeta\in\dom(H^{k+1})$. This yields $x_1\in\dom(H^{k+1})+\rrm$ since $\zeta\in\rrm$. Note that $x=x_1+x_2$ where $x_2\in\rrm$. As such, we can conclude that $x\in \big(\dom(H^{k+1})+\rrm\big)\cap\calV$. Finally, we take $\calV=\rrp$ in Lemma~\ref{lemm:dom he dom h} to conclude that $x\in\big(\dom(H^{k+1})+\rrp\big)\cap\calV$ which proves that $\dom(\Ho^{k+1})\subseteq \big(\dom(H^{k+1})+\rrp\big)\cap\calV$.

{\ref{lemmitem:feasible H}}: From the statement \ref{lemmitem:ho=L+out}, we know that $\Ho$ is a linear process. Therefore, we have
\beq\label{e:char hout domain}
\calF(\Ho)=\dom(\Ho^n)
\eeq
due to Lemma~\ref{lemm:linear F and R}. Clearly, $\calF(H)\subseteq\dom(H^n)$. The reverse inclusion would follow from Lemma~\ref{lemm:Ffdiff} if $\dom(H^n)\subseteq \dom (H^{n+1})$. 

Let $\xi\in\dom(H^n)$. As such, we see that there exist $x_0,x_1,\ldots,x_n$ such that $\xi=x_0$ and $x_{k+1}\in H(x_k)$ for all $k$ with $0\leq k \leq n-1$. From the direct sum \eqref{eq:direct sum R(L+) and W}, there exist $y_0,y_1,\ldots,y_n\in\rrp$ and $z_0,z_1,\ldots,z_n\in\calV$ such that $x_k=y_k+z_k$ for all $0\leq k \leq n-1$. Note that $(z_k,z_{k+1})=(x_k,x_{k+1})-(y_k,y_{k+1})\in\gr(\Ho)$ since $(x_k,x_{k+1})\in\gr(H)$ and $(y_k,y_{k+1})\in\rrp\times\rrp$. 
 
As $\Ho$ is single valued, and $z_0\in\calF(\Ho)$, we know that after $n$ steps, we are still inside the feasible set: $z_n\in\calF(\Ho)\subseteq\dom(\Ho)$. The domain condition and $\dom(H)+\rrm=\dom(H)+\rrp$ and \ref{lemmitem:domHout} imply that $\dom(\Ho)=\big(\dom( H )+\rrm\big)\cap \calV$. Therefore there exist $\bar{z}\in\dom(H)$ and $\tilde{z}\in\rrm$ such that $z_n=\barz+\tilde{z}$. Due to the domain condition and Lemma~\ref{lemm:conespacesequiv}.\ref{lemmitem:subplus} we have $\rrp\subseteq\dom(H)+\rrm$. Then, there exist $\bar{y}\in\dom(H)$ and $\tilde{y}\in\rrm$ such that $y_n=\bary+\tilde{y}$. We know that $\rrm=L^n_-(0)\subseteq H^n(0)$ where the last inclusion follows from \eqref{eq:lh-+}. This means that $\bar{y}+\bar{z}=x_n-\tilde{y}-\tilde{z}\in H^n(x_0)+H^n(0)\subseteq H^{n}(x_0)$. Since $\bar{y}+\bar{z}\in\dom(H)$, we see that $\xi=x_0\in\dom(H^{n+1})$. Therefore, we see that $\dom(H^n)\subseteq \dom(H^{n+1})$, which proves the statement.

{\ref{lemmitem:F hout}}: From \eqref{e:char hout domain}, \ref{lemmitem:ho^k} and \ref{lemmitem:feasible H}, we see that $\calF(\Ho)=\dom(\Ho^n)=\big(\calF(H)+\rrp\big)\cap\calV$. Since $\Ho$ is a linear process due to \ref{lemmitem:ho=L+out}, $\calF(\Ho)$ is a subspace.

{\ref{lemmitem: calfh + rrp is subspace}}: Since $\calF(H)+\rrp$ is a convex cone, it is enough to show that $\xi\in\calF(H)+\rrp$ implies $-\xi\in\calF(H)+\rrp$. Let $\xi\in\calF(H)+\rrp$. Then, $\xi=\zeta+\eta$ where $\zeta\in \calF(H)$ and $\eta\in\rrp$. We also know from \eqref{eq:direct sum R(L+) and W} that $\xi=\xi_1+\xi_2$ where $\xi_1\in\rrp$ and $\xi_2\in\calV$. It now follows from \ref{lemmitem:F hout} that $\xi_2=\zeta+\eta-\xi_1\in\calF(\Ho)$. Since $\calF(\Ho)$ is a subspace, we have that $-\xi_2=\calF(H)+\rrp$. Then, $-\xi=-\xi_1-\xi_2\in\calF(H)+\rrp$. Therefore, $\calF(H)+\rrp$ is a subspace.
\EP

It turns out that the trajectories of the difference inclusion \eqref{eq:incl} can be decomposed according to \eqref{eq:direct sum R(L+) and W} by using the outer process $\Ho$ as stated next.

\blem\label{l:Kalman-like}
Let $H$ be a convex process and $(x_k)_{k\in\N}\in\mathfrak{B}(H)$. Then, there exist sequences $(\xi_k)_{k\in\N}\subset\Fo$ and $(\eta_k)_{k\in\N}\subset\rrp$ such that $x_k=\xi_k+\eta_k$ and $\xi_{k+1}=\res{\Ho}{\Fo}(\xi_k)$ for all $k\in\N$. 
\elem

\BP
Existence of sequences $(\xi_k)_{k\in\N}\subset\calV$ and $(\eta_k)_{k\in\N}\subset\rrp$ such that $x_k=\xi_k+\eta_k$ for all $k\in\N$ follows from the direct sum \eqref{eq:direct sum R(L+) and W}. From \eqref{eq:def H out}, we see that $\xi_{k+1}\in\Ho(\xi_k)$ and hence $\xi_k\in\Fo$ for all $k\in\N$. Since $\Ho$ is a linear map due to Lemma~\ref{lemm:Houtthings}.\ref{lemmitem:Loutsinglevalued}-\ref{lemmitem:ho=L+out}, we further see that $\xi_{k+1}=\res{\Ho}{\Fo}(\xi_k)$ for all $k\in\N$.
\EP

\section{Proofs}\label{sec:proofs}

This section will use the derived framework, and previously proven results to prove the main results of this paper.

\subsection{Proof of Lemma~\ref{l:ho single valued lin proc}}
It follows from Lemma~\ref{lemm:Houtthings} that $\Ho$ is a single valued linear process, i.e. a linear map and $\calF(\Ho)=\big(\calF(H)+\rrp\big)\cap\calV$ is a subspace. Since $\calF(\Ho)$ is weakly $\Ho$ invariant due Lemma~\ref{lemm:feasible largest wHi cone}, single valuedness of $\Ho$ readily implies that $\Ho(\Fo)\subseteq\Fo$.

\subsection{Proof of Theorem~\ref{thm:main-r}}
We will prove the implications \ref{thm:main-r1} $\Rightarrow$ \ref{thm:main-r2}, \ref{thm:main-r2} $\Rightarrow$ \ref{thm:main-r3}, and \ref{thm:main-r3} $\Rightarrow$ \ref{thm:main-r1}.

{\em \ref{thm:main-r1} $\Rightarrow$ \ref{thm:main-r2}:} From \eqref{e:f hi}, we have that $\calF(\Hi)=\calF(H)\cap\rrp$. Since $H$ is reachable, this implies that $\calF(\Hi)\subseteq\calR(H)$. In view of \eqref{e:r hi= r h}, we have that $\calF(\Hi)\subseteq\calR(\Hi)$. Therefore, $\Hi$ is reachable. To show that $\Ho$ is reachable, note that $\calF(\Ho)=\big(\calF(H)+\rrp\big)\cap\calV$ due to Lemma~\ref{lemm:Houtthings}.\ref{lemmitem:F hout}. Since $H$ is reachable, we see that $\calF(\Ho)\subseteq\rrp\cap\calV=\zset$. As such, $\Ho$ is reachable.  

{\em \ref{thm:main-r2} $\Rightarrow$ \ref{thm:main-r3}:} From Lemma~\ref{l:Hin char}.\ref{l:Hin char.r}, reachability of $\Hi$ implies that all eigenvectors of $\Hi^-$ corresponding to eigenvalues in $[0,\infty)$ belong to $\rrp^\bot$. Note that $\Ho$ is single valued due to Lemma~\ref{lemm:Houtthings}.\ref{lemmitem:Loutsinglevalued} and \ref{lemmitem:ho=L+out}. This implies that $\Ho(0)=\{0\}$ and hence that $\calR(\Ho)=\zset$. Therefore, reachability of $\Ho$ implies that $\calF(\Ho)=\zset$. 

{\em \ref{thm:main-r3} $\Rightarrow$ \ref{thm:main-r1}:} Since $\calF(\Ho)=\zset$, 
it follows from \eqref{eq:direct sum R(L+) and W} and Lemma~\ref{lemm:Houtthings}.\ref{lemmitem:F hout} that $\calF(H)\subseteq\rrp$. Then, we have that $\calF(H)=\calF(\Hi)$ in view of \eqref{e:f hi}. From Lemma~\ref{l:Hin char}.\ref{l:Hin char.r}, we see that $\Hi$ is reachable since all eigenvectors of $\Hi^-$ corresponding to eigenvalues in $[0,\infty)$ belong to $\rrp^\bot$. Therefore, we have $\calF(H)=\calF(\Hi)\subseteq\calR(\Hi)$. Since $\calR(\Hi)=\calR(H)$ due to \eqref{e:r hi= r h}, $H$ is reachable.

For the last claim, suppose that $H$ is reachable. From Theorem~\ref{thm:main-r}.\ref{thm:main-r2}, we know that $\Hi$ is reachable. In turn, Lemma~\ref{l:Hin and W} (implication \ref{l:Hin and W.1} $\Rightarrow$ \ref{l:Hin and W.2} for $\calW=\zset$) implies that $\calR(\Hi)=\rrp$. Since it readily holds that $\calR(\Hi)\subseteq\calR(H)\subseteq\rrp$, we  see that $\calR(H)=\rrp$. By taking $\calC_\ell=H^\ell(0)$ and applying Lemma~\ref{lemm:FGisFD}, we obtain $\calR(H)= H^q(0)$ for some $q\geq 0$.

\subsection{Proof of Theorem~\ref{thm:main-s}} The implication {\em \ref{thm:main-s2} $\Rightarrow$ \ref{thm:main-s1}} is evident. In what follows, we will prove the implications {\em \ref{thm:main-s1} $\Rightarrow$ \ref{thm:main-s4}}, {\em \ref{thm:main-s4} $\Rightarrow$ \ref{thm:main-s3}}, and {\em \ref{thm:main-s3} $\Rightarrow$ \ref{thm:main-s2}}.

{\em \ref{thm:main-s1} $\Rightarrow$ \ref{thm:main-s4}}: Since $H$ is stabilizable, it follows from Lemma~\ref{l:Kalman-like} that all eigenvalues of the linear map $\res{\Ho}{\Fo}$ are in the open unit disc. To prove the rest, we first observe that $\Hi$ is stabilizable whenever so is $H$. Now, let $\lambda\geq 1$ and $\xi$ be such that $\lambda\xi\in\Hi^-(\xi)$. The arguments used in the proof of Lemma~\ref{l:Hin char}.\ref{l:Hin char.s} result in $\xi\in\big(\calR(\Hi)-\mathcal{S}(\Hi)\big)^-$. From Lemma~\ref{l:invariance s}, we know that $\Hi\inv(\mathcal{S}(\Hi))\subseteq\mathcal{S}(\Hi)$. By applying Lemma~\ref{l:Hin and W} (implication \ref{l:Hin and W.4} $\Rightarrow$ \ref{l:Hin and W.2} for $\calW=\mathcal{S}(\Hi)$), we see that $\calR(\Hi)-\mathcal{S}(\Hi)=\rrp$. Therefore, we see that $\xi\in\rrp^\bot$.

{\em \ref{thm:main-s4} $\Rightarrow$ \ref{thm:main-s3}}: We can conclude that $\Hi$ is exponentially stabilizable by applying Lemma~\ref{l:Hin char}.\ref{l:Hin char.s}. For the second part, we know by Lemma~\ref{l:ho single valued lin proc} that $\Ho$ is a single valued linear map on $\Fo$. As such, it is exponentially stabilizable if and only if all eigenvalues of that linear map are in the open unit disc.

{\em \ref{thm:main-s3} $\Rightarrow$ \ref{thm:main-s2}}: To prove this implication, we will first construct a bounded polyhedron $\calP$ as in Lemma~\ref{lemm:poly to stab} and then show that every feasible point can be steered to $\calP$ in finitely many steps. 

To construct $\calP$, we need some preparation. Let $\calW$ be a subspace such that 
\beq\label{e:decomp rrp}
\calR_+ = (\Lin(\calF(H))\cap \calR_+)\oplus \calW.
\eeq
For $\xi\in\calF(\Ho)$, define
\beq\label{e:def G}
G(\xi) := \Big( \big( (\calF(H)-\xi)\cap \calR_+\big)+\big( \Lin(\calF(H))\cap \calR_+\big) \Big) \cap \calW.
\eeq
It can be easily verified that $G$ is a convex process. From Lemma~\ref{lemm:Houtthings}.\ref{lemmitem:F hout}, we know that for every $\xi\in\calF(\Ho)$ there exists $\eta\in\calR_+$ such that $\xi+\eta\in \calF(H)$. Further, the decomposition \eqref{e:decomp rrp} implies that there exist $\eta_1\in  \Lin(\calF(H))\cap \calR_+$ and $\eta_2\in\calW$ such that $\eta=\eta_1+\eta_2$. Then, we see that $\eta_2\in G(\xi)$. Therefore, we have that $\dom G =\calF(\Ho)$. 

Now, we claim that $G$ is single valued. For this, take $\eta,\zeta\in G(\xi)$. This means that there exist $\eta_1,\zeta_1\in\calR_+$ and $\eta_2,\zeta_2\in \Lin(\calF(H))\cap \calR_+$ such that
\[ \xi+\eta_1\in \calF(H), \quad \xi+\zeta_1\in \calF(H),\quad \eta_1=\eta_2+\eta, \quad \zeta_1=\zeta_2+\zeta.\]
By using the first two relations, we obtain $\eta_1-\zeta_1 \in \Lin(\calF(H))$. This implies that $\eta_1-\zeta_1 \in \Lin(\calF(H))\cap \calR_+$. Note that $\eta-\zeta=(\eta_1-\zeta_1)-(\eta_2-\zeta_2)$. Therefore, we see that $\eta-\zeta\in \Lin(\calF(H))\cap \calR_+$. Since  $\eta-\zeta\in \calW$, it follows from \eqref{e:decomp rrp} that $\eta=\zeta$, i.e. $G$ is single valued. Since its domain is a subspace, we further see that $G$ is a linear map.

Since $\Fo$ is a subspace due to Lemma~\ref{lemm:Houtthings}.\ref{lemmitem:Loutsinglevalued}-\ref{lemmitem:ho=L+out}, we can find a bounded polyhedron $\calQ\subseteq\Fo$ containing the unit ball relative to $\Fo$, i.e. 
\beq\label{e:unit ball in calQ}
\set{\xi\in\Fo}{\abs{\xi}\leq 1}\subset\calQ.
\eeq
From \eqref{e:def G}, we see that 
\beq\label{e:I+G calQ in Lin FH}
(I+G)\calQ\subseteq \Lin(\calF(H)).
\eeq
Moreover, $(I+G)\calQ$ is a bounded polyhedron since $G$ is a linear map. 

From Lemma~\ref{lemm:RH invariance}, we know that $\calR(H)-\calR(H)=\rrp$. Then, Lemma~\ref{lemm:FGisFD} implies that there exists $q_1\geq 0$ such that 
\beq\label{e: h0-h0}
H^q(0)-H^q(0)=\rrp
\eeq
for all $q\geq q_1$. Since $\Hi$ is exponentially stabilizable, we have 
$\calS_{\mathrm{e}}(\Hi) -\calR(H)=\rrp$ due to 
Lemma~\ref{l:Hin and W} (implication \ref{l:Hin and W.4} $\Rightarrow$ \ref{l:Hin and W.2} for $\calW=\mathcal{S}_{\mathrm{e}}(\Hi)$) and the fact that $\calR(\Hi)=\calR(H)$. By applying Lemma~\ref{lemm:FGisFD}, we see that there exists $q_2\geq 0$ such that $\calS_{\mathrm{e}}(\Hi) -H^q(0)=\rrp$ for $q\geq q_2$. 

Let $q_3=\max\{q_1,q_2\}$. Since $\calS_{\mathrm{e}}(\Hi) -H^{q_3}(0)$ is a subspace, it follows from Lemma~\ref{l:cone min cone is subspace} that there exists $y$ such that
\beq\label{e:y in ri se-h0}
y\in \rint (\calS_{\mathrm{e}}(\Hi))\cap \rint(H^{q_3}(0)).
\eeq
As $\rint (\cone(y))\cap \rint(H^{q_3}(0))\neq \emptyset$, we know $\cone(y)-H^{q_3}(0)$ is a subspace from Lemma~\ref{l:cone min cone is subspace}. As such, it contains both $-H^{q_3}(0)$ and $H^{q_3}(0)$. This means that $\calR_+ = \Lin (H^{q_3}(0)) \subseteq \cone(y)-H^{q_3}(0)\subseteq \calR_+$. Thus, we see that
\beq\label{e:cone y - h q}
\cone(y)-H^{q_3}(0)=\rrp.
\eeq
As $H^{q_3}(0)\subseteq H^q(0)$ for all $q\geq q_3$ and $q_3\geq q_1$, \eqref{e: h0-h0} implies that 
\beq\label{e:ri h q grows}
\rint(H^{q_3}(0))\subseteq \rint(H^q(0))
\eeq
for all $q\geq q_3$. Therefore, we see that
\begin{gather}
y\in \rint(\calS_{\mathrm{e}}(\Hi))\cap \rint(H^q(0)),\label{e: y in ri q3}\\
\cone(y)-H^q(0)=\rrp\label{e: cone y - h 0=rrp}
\end{gather} 
for all $q\geq q_3$. 

Now, let $\lambda\in(0,1)$. Let $\calB=\{x\in\calR_+ \mid \abs{x}\leq 1\}$ denote the unit ball in $\calR_+$. From \eqref{e:ri h q grows} and \eqref{e: y in ri q3}, we see that there exists $\varepsilon>0$ such that $y+\varepsilon \calB\subseteq H^q(0)$ for all $q\geq q_3$. Since $y\in \calS_{\mathrm{e}}(\Hi)$, 
there exists $q_4\geq 0$ such that for all $q\geq q_4$ there exists $y_q\in H^q(y)$ with $\abs{y_q}\leq \varepsilon\lambda$. As $y\in\rrp$ and $\rrp$ is strongly $H$ invariant, $y_q\in\rrp$. Then, $\lambda y-y_q \in H^q(0)$ for all $q\geq q_4$. As such, we can conclude that 
\beq\label{e: lambda is eig for H with y}
\lambda y = y_q-y_q+\lambda y \in H^q(y)+H^q(0) \subseteq H^q(y) \eeq
for all $q\geq q_4$.

Since $\Hi$ is exponentially stabilizable, $\calF(\Hi)=\calS_{\mathrm{e}}(\Hi))$. In view of \eqref{e: y in ri q3}, this means that $y\in\rint(\calF(\Hi))$. From Lemma~\ref{lemm:Houtthings}.\ref{lemmitem: calfh + rrp is subspace}, we have that $\calF(H) + \calR_+$ is a subspace. According to Lemma~\ref{l:cone min cone is subspace}, we have that $\rint(\calF(H))\cap \calR_+\neq \emptyset$. Together with \cite[Thm. 6.5]{Rockafellar:70},  this implies that $\rint(\calF(H)\cap \calR_+) = \rint(\calF(H))\cap \calR_+$. Thus, we see that 
\[y\in\rint(\calF(\Hi))= \rint(\calF(H)\cap \calR_+) = \rint(\calF(H))\cap \calR_+.\]
Now, let $\calQ= \conv \{\xi_0^i \mid i = 1,\ldots, r\}$. Since $y\in\rint(\calF(H))$, it follows from \eqref{e:I+G calQ in Lin FH} that there exists $\beta_1>0$ such that 
\beq\label{e:xi i 0 brought in fh}
(I+G)\xi_0^i+\beta y\in\calF(H)
\eeq
for all $i$ and $\beta\geq\beta_1$. Let $\xi^i_q\in \Ho^q(\xi_0^i)$. Take $\rho\in(\lambda,1)$. As $\Ho$ is single valued and exponentially stabilizable, there exists $q_5\geq 0$ such that $\abs{\xi^i_q}<\rho$ for all $i$ and $q\geq q_5$. 

From \eqref{e:xi i 0 brought in fh}, we see that there exists $(z^i_k)_{k\in\N}\in\mathfrak{B}(H)$ with $z^i_0=\xi^i_0 +G(\xi^i_0)+\beta y$. In view of \eqref{eq:direct sum R(L+) and W} and \eqref{e:decomp rrp}, we have that 
\[\mathbb{R}^n = (\Lin(\calF(H))\cap \calR_+)\oplus \calW \oplus \calV. \]
To this, we can apply Lemma~\ref{l:Kalman-like} and the definition of $G$ in \eqref{e:def G}, to note that
\[ z^i_k = \xi^i_k +G(\xi^i_k) +\eta^i_k, \]
where $\xi^i_k\in\Ho^k(\xi^i_0)$ and $\eta^i_k\in\rrp$ for all $k\in\N$. Now, let $q=\max\pset{q_3,q_4,q_5}$. Then, we see from \eqref{e: cone y - h 0=rrp} that $\eta^i_q = \alpha^iy-\zeta^i$, where $\alpha^i\geq 0$ and $\zeta^i\in H^q(0)$. Then, we can apply \eqref{eq: H(x) + H(0)= H(x)} to obtain that
\[ \xi^i_q +G(\xi^i_q) + \alpha^iy \in H^q\big(\xi^i_0+G(\xi_0^i)+ \beta y\big). \] 
From \eqref{e: y in ri q3} and \eqref{e: lambda is eig for H with y} we know that $y\in H^q(0)$ and $\lambda y\in H^q(y)$. Therefore, we have
\begin{equation}\label{eq:construction}
\xi^i_q +G(\xi^i_q) + (\alpha^i+\gamma\lambda+\delta^i)y \in H^q\big(\xi^i_0+G(\xi_0^i)+ (\beta+\gamma) y\big)
\end{equation}
for all $\gamma\geq 0$ and $\delta^i\geq 0$. Note that we can take $\gamma$ and $\delta^i$ such that $\alpha^i+\gamma\lambda+\delta^i=\rho (\beta+\gamma)$ for all $i$ since $\lambda<\rho$. Then,  \eqref{eq:construction} boils down to 
\[ \xi^i_q +G(\xi^i_q) + \rho (\beta+\gamma)y \in H^q\big(\xi^i_0+G(\xi_0^i)+ (\beta+\gamma) y\big). \] 

Now, take $\calP= (I+G)\calQ+\bar{\beta} y$ where $\bar{\beta}= \beta+\gamma$. Let $x\in \calP$. Then, $x=\xi_0+G(\xi_0)+\bar{\beta}y$ where $\xi_0\in\calQ$. This means that $\xi_0 = \sum_i{a_i\xi^i_0}$ where $a_i\geq 0$ and $\sum_i{a_i} =1$. Thus, we can conclude that 
\[ \xi+G(\xi) + \rho \bar{\beta}y \in H^q\big(\xi_0+G(\xi_0)+\bar{\beta} y\big)= H^q(x) \]
where $\xi= \sum_i{a_i\xi^i_q}$. Since $\abs{\xi^i_q}\leq \rho$, we have $\abs{\xi}\leq \rho$. As such, we see from \eqref{e:unit ball in calQ} that $\xi\in \rho \calQ$. This proves that $H^q(x)\cap \rho \calP \neq \emptyset$. By applying Lemma~\ref{lemm:poly to stab}, we can conclude that $\calP\subseteq\calS_{\mathrm{e}}(H)$.

Now, let $\barx\in\calF(H)$. If $\barx=0$, then we clearly have that $\barx\in\calS_{\mathrm{e}}(H)$. Suppose that $\barx\neq0$. Let $(x_k)_{k\in\N}\in\mathfrak{B}(H)$ be a trajectory with $x_0=\barx$. From \eqref{eq:direct sum R(L+) and W}, \eqref{e:decomp rrp}, and Lemma~\ref{l:Kalman-like}, we see that $x_k=\xi_k+G(\xi_k)+\eta_k$ where $\xi_{k+1}\in\res{\Ho}{\Fo}(\xi_k)$ and $(\eta_k)_{k\in\N}\subset\rrp$. Since $\Ho$ is exponentially stabilizable, there exists $\barq\geq\max\pset{q_3,q_4,q_5}$ such that $\abs{\xi_{\barq}}\leq 1$. From \eqref{e: y in ri q3}, \eqref{e: cone y - h 0=rrp} and \eqref{eq: H(x) + H(0)= H(x)} we see that there exists $\bar{\alpha}$ such that 
\beq\label{e:xi barq}
\xi_{\barq}+G(\xi_{\barq})+\alpha y\in H^{\barq}(\barx)
\eeq
for all $\alpha\geq\bar{\alpha}$. Due to \eqref{e:unit ball in calQ}, we have that
$\xi_{\barq}+G(\xi_{\barq})\in(I+G)\calQ$.
As such, $\xi_{\barq}+G(\xi_{\barq})+\bar{\beta}y\in\calP\subseteq\calS_{\mathrm{e}}(H)$. By taking $\alpha'=\max(\bar{\alpha},\bar{\beta})$, we see that $\xi_{\barq}+G(\xi_{\barq})+\alpha' y\in\calS_{\mathrm{e}}(H)$ since $y\in\calS_{\mathrm{e}}(H)$ due to \eqref{e: y in ri q3}. From \eqref{e:xi barq}, we have that 
$\barx\in H^{-\barq}(\xi_{\barq}+G(\xi_{\barq})+\alpha' y)$. Therefore, it follows from Lemma~\ref{l:invariance s} that $\barx\in\calS_{\mathrm{e}}(H)$. This proves that $\calF(H)\subseteq\calS_{\mathrm{e}}(H)$ and hence $H$ is exponentially stabilizable.\EP

\subsection{Proof of Theorem~\ref{t:reach imply null and cont}}
Clearly, $H$ is controllable if and only if $H$ is both reachable and null-controllable. Therefore, what needs to be proven is that reachability of $H$ implies its null-controllability. From Theorem~\ref{thm:main-r} we obtain that $\calR(H)=\rrp$ and $\calR(H)= H^q(0)$ for some $q\geq 0$. Let $\xi\in\calF(H)\subseteq \calR(H)$. Then, there exists a trajectory $(x_k)_{k\in\N}\in\calB(H)$ with $x_0=\xi$. Clearly, we have $x_q\in H^q(\xi)$. As $\calR(H)$ is strongly $H$ invariant, we know that $x_q\in\calR(H)=\rrp$. Hence, we see that $-x_q\in\calR(H)=H^q(0)$. It then follows from \eqref{eq: H(x) + H(0)= H(x)} that $0= x_q-x_q \in H^q(\xi)+H^q(0)=H^q(\xi)$. Consequently, $\calF(H)\subseteq\calN(H)$, that is $H$ is null-controllable.

\subsection{Proof of Theorem~\ref{thm:main-n}} We will prove the implications {\em \ref{thm:main-n1} $\Rightarrow$ \ref{thm:main-n3}}, {\em \ref{thm:main-n3} $\Rightarrow$ \ref{thm:main-n2}}, and {\em \ref{thm:main-n2}} $\Rightarrow$ \ref{thm:main-n1}.

{\em \ref{thm:main-n1} $\Rightarrow$ \ref{thm:main-n3}}: Since $H$ is null-controllable, it follows from Lemma~\ref{l:Kalman-like} that all trajectories of $\res{\Ho}{\Fo}$ reach the origin in a finite number of steps. As such, all eigenvalues of the linear map $\res{\Ho}{\Fo}$ are zero and thus $\res{\Ho}{\Fo}$ is nilpotent. For the rest, we first observe that $\Hi$ is null-controllable whenever so is $H$ due to 
\eqref{e:f hi} and \eqref{e:n hi}. Then, it follows from Lemma~\ref{l:Hin char}.\ref{l:Hin char.n1} that all eigenvectors of $\Hi^-$ corresponding to eigenvalues in $(0,\infty)$ belong to $\rrp^\bot$.

{\em \ref{thm:main-n3} $\Rightarrow$ \ref{thm:main-n2}}: This implication follows from Lemma~\ref{l:Kalman-like} and Lemma~\ref{l:Hin char}.\ref{l:Hin char.n2} in the same manner as the proof of Theorem~\ref{thm:main-s} ({\em \ref{thm:main-s4} $\Rightarrow$ \ref{thm:main-s3}}).

{\em \ref{thm:main-n2} $\Rightarrow$ \ref{thm:main-n1}}: Let $\barx\in\calF(H)$. Then, there exists a trajectory $(x_k)_{k\in\N}\in\mathfrak{B}(H)$ such that $x_0=\barx$. From Lemma~\ref{l:Kalman-like}, we have that $x_k=x_k=\xi_k+\eta_k$ where $\eta_k\in\rrp$ and $\xi_{k+1}\in\res{\Ho}{\Fo}(\xi_k)$ for all $k\in\N$. Since $\Ho$ is null-controllable, there must exist $q\geq 0$ such that $\xi_q=0$. This means that $\eta_q\in\calF(H)\cap\rrp=\calF(\Hi)$ since $x_q\in\calF(H)$. As such, we have that $x_q\in\calN(\Hi)\subseteq\calN(H)$ since $\Hi$ is null-controllable. This means that $\barx\in H^{-q}(x_q)\subseteq H^{-q}(\calN(H))$. Then, it follows from Lemma~\ref{l:invariance n} that $\barx\in\calN(H)$. This proves that $\calF(H)\subseteq\calN(H)$ and hence $H$ is null-controllable.

\section{Conclusion}\label{sec:conclusion}
In this paper, we have developed a framework for analysis of convex processes. Central concepts in this are weakly and strongly invariant cones, the minimal and maximal linear processes and duality. It was shown that these concepts naturally have a central role in the analysis of convex processes. 

Within this, we developed Hautus-type spectral tests for reachability, stabilizability and null-controllability of nonstrict convex processes. In essence we have shown that, under a condition on the domain, we can investigate the properties of a convex process by considering its \textit{inner} and \textit{outer} processes separately. This result is akin to the so-called Kalman decomposition for linear systems. After this, spectral characterizations for either of these processes were developed. In particular for the inner process, this required additional developments in duality of convex processes. Moreover, we have proven that, under the domain condition, reachability and controllability are equivalent. 

It was shown that these main results generalize all previously known characterizations. In particular, the known results for (strict) convex processes and for linear processes were unified.

\subsection*{Future work}
As noted, the framework of this paper can be applied to many different problems regarding convex processes. The results in duality will prove useful for different stabilizability problems. Indeed, it was shown in \cite{Eising2020b}, that a stricter domain condition plays a role in the study of (duality of) Lyapunov functions for convex processes. 

Another avenue of extensions would be to consider control of convex processes. While in the case of linear systems stabilizability is equivalent to the existence of a static state feedback controller, for convex processes such results do not hold. An interesting problem is to develop a theory of control for this class of systems. 

In this paper, we considered convex processes, with a motivation of developing results for conically constrained linear systems. A logical extension would be to investigate properties of more general set-valued maps. In particular set-valued maps with a convex graph would prove interesting. Parallel to the work of this paper, an investigation of properties of set-valued maps with affine graphs would prove another relevant intermediate point towards the general convex case. The first steps in this regard were made in \cite{Kaba:18}.

\bibliography{convex_processes}
\bibliographystyle{siam}

\end{document}